%% file: merged_paper.tex
\documentclass[journal]{IEEEtran}
\IEEEoverridecommandlockouts

\usepackage{graphicx}
\usepackage{cite}
\usepackage{amsfonts} 
\usepackage{amsmath}  
\usepackage{bm}
\usepackage{psfrag}
\usepackage{amsbsy}
\usepackage{amssymb}
\usepackage{mathrsfs}
\usepackage{stfloats}
\usepackage{dsfont}
\usepackage{setspace}
\usepackage{array}
\usepackage{bbm}
\usepackage{mathabx}
\usepackage{algorithmic}
\usepackage{algorithm}
\usepackage{glossaries}
\usepackage{color}
\usepackage{amsthm}
\usepackage{lettrine}

\newtheoremstyle{mystyle}
{}                
{}                
{\itshape}        
{}                
{\bfseries}       
{:}               
{.5em}            
{}                

\theoremstyle{mystyle}
\newtheorem{proposition}{Proposition}

\newtheorem{remark}{Remark}
\newcommand{\ds}{\displaystyle}

\DeclareMathOperator{\tr}{tr}

\def\BibTeX{{\rm B\kern-.05em{\sc i\kern-.025em b}\kern-.08em
    T\kern-.1667em\lower.7ex\hbox{E}\kern-.125emX}}

\input{style_icc2024_ext.tex}
\input{acronyms}

\title{Secrecy Energy Efficiency Maximization in RIS-Aided Networks: Active or Nearly-Passive RIS?}
\author{
	Robert Kuku Fotock, {\em Student Member, IEEE}, Agbotiname Lucky Imoize, {\em Senior Member, IEEE},  Alessio Zappone,\\ {\em Fellow, IEEE}, Marco Di Renzo, {\em Fellow, IEEE}, Roberto Garello, {\em Senior Member, IEEE}
	\thanks{R. K. Fotock and A. Zappone are with CNIT and with University of Cassino and Southern Lazio, Italy (robert.fotock@cnit.it; alessio.zappone@unicas.it). A. L. Imoize is with CNIT and with Politecnico di Torino, Italy (aimoize@cnit.it). M. Di Renzo is with Universit\'e Paris-Saclay, CNRS, CentraleSup\'elec, Laboratoire des Signaux et Syst\`emes, 3 Rue Joliot-Curie, 91192 Gif-sur-Yvette, France. (marco.di-renzo@universite-paris-saclay.fr), and with King's College London, Centre for Telecommunications Research -- Department of Engineering, WC2R 2LS London, United Kingdom (marco.di\_renzo@kcl.ac.uk).. R. Garello is with Politecnico di Torino, Italy (roberto.garello@polito.it).
		The work of R. K. Fotock was supported by the European Commission through the H2020-MSCA-ITN-METAWIRELESS project until November 30th, 2024. Afterwards, the work of R. K. Fotock was supported by the Project “GARDEN”, funded by EU in NextGeneration EU plan, Mission 4, Component 1, CUP H53D23000480001, through the Italian “Bando Prin 2022 - D.D. 104 del 02-02-2022“ by MUR. The work of A. L. Imoize was supported by the European Commission through the HE-MSCA-DN-INTEGRATE project, grant agreement 101072924. The work of A. Zappone has been funded by the European Union - NextGenerationEU under the project NRRP RESTART, RESearch and innovation on future Telecommunications systems and networks, to make Italy more smART PE\_00000001 - Cascade Call SPARKS project, with CUP D43C22003080001. The work of M. Di Renzo was supported in part by the European Union through the Horizon Europe project COVER under grant agreement number 101086228, the Horizon Europe project UNITE under grant agreement number 101129618, the Horizon Europe project INSTINCT under grant agreement number 101139161, and the Horizon Europe project TWIN6G under grant agreement number 101182794, as well as by the Agence Nationale de la Recherche (ANR) through the France 2030 project ANR-PEPR Networks of the Future under grant agreement NF-PERSEUS 22-PEFT-004, and by the CHIST-ERA project PASSIONATE under grant agreements CHIST-ERA-22-WAI-04 and ANR-23-CHR4-0003-01. This work was partially presented in \cite{FotockSPAWC24}. 
}}

\begin{document}

\maketitle

\begin{abstract}
This work addresses the problem of secrecy energy efficiency (SEE) maximization in RIS-aided wireless networks. The use of active and nearly-passive RISs are compared and their trade-off in terms of SEE is analyzed. Considering both perfect and statistical channel state information, two SEE maximization algorithms are developed to optimize 
 the transmit powers of the mobile users, the RIS reflection coefficients, and the base station receive filters. Numerical results quantify the trade-off between active and nearly-passive RISs in terms of SEE, with active RISs yielding worse SEE values as the static power consumed by each reflecting element increases. 
\end{abstract}

\vspace{1em}

\begin{IEEEkeywords}
RIS, secrecy energy efficiency, physical layer security, resource allocation, perfect and statistical CSI.
\end{IEEEkeywords}

\IEEEpeerreviewmaketitle

\section{Introduction} 
\lettrine[nindent=0.1em,lines=2]{R}{ECONFIGURABLE} intelligent surfaces (RISs) have emerged as key enablers of future 6G wireless networks ~\cite{di2020analytical, SmartWireless, RuiZhang_COMMAG, huang2019indoor,YuanStandard}. RISs have two major advantages over other available technologies: (1) they require simpler hardware, which allows the use of a large number of electromagnetic elements, in turn leading to higher rates; (2) they require a lower energy expenditure, fully operating in the analog domain and being, at least in their first version, nearly-passive devices. As a result, RISs can provide unprecedented levels of energy efficiency (EE), which is a major requirement of 6G wireless networks. Indeed, EE already was a major requirement of 5G and is even more relevant for future 6G networks, especially considering that 5G networks did not achieve the 2000x EE increase that was sought after \cite{David5G}. In fact, while 5G technologies like massive multiple-input multiple-output (MIMO) provide much higher rate levels, it is argued in \cite{Huawei2020} that they also consume up to three times more energy than legacy 4G  technologies, mainly due to the large size of digital antenna arrays. As a result, the EE of 3rd Generation Partnership Project (3GPP) new radio deployments is estimated to be approximately four times larger than that of 3GPP long term evolution deployments.

On the other hand, RISs suffer from the so-called double or multiplicative fading effect, i.e. the end-to-end channel from the transmitter, through the RIS, to the receiver, is the product of the channels of the transmitter-RIS and RIS-receiver links. In order to compensate for this effect and obtain satisfactory rate performance, the use of an exceedingly large number of reflecting elements may be required \cite{Dajer2022,Dong2022,Sihlbom2022}, which is not ideal as far as both energy consumption and computational complexity are concerned. For this reason, the use of active RISs has been proposed \cite{Long2021}. An active RIS consists of reflectors characterized by active-load impedances, which can reflect and amplify the incident signal impinging on the RIS. Unlike traditional relays, the amplification occurs at the electromagnetic level and thus requires neither a radio frequency chain nor any conversion between the digital and analog domains. Active RISs can combat the multiplicative fading effect \cite{Zhang2023}. However, while active RISs can provide higher rates compared to their nearly-passive counterparts, their performance in terms of network EE is not clear, because the use of a power amplifier is an additional source of power consumption that tends to reduce the energy efficiency of the network. Moreover, active RISs also amplify thermal noise due to the presence of power amplifiers. 

Together with the EE, another key performance requirement of future 6G networks is the  privacy of the data transfer \cite{Ziegler2021,Norrman2024}. Due to its inherent nature, the wireless channel is prone to being overheard by other nodes of the network and wiretapped by malicious users. Among the different techniques for the security of wireless networks, physical layer security has the advantage of not requiring any dedicated cryptographic algorithm, relying instead on information-theoretic security, thus reducing complexity and communication delays. Physical layer security can be used in conjunction with traditional cryptography to provide an additional layer of security \cite{Mitev2022}, and is considered a promising way of providing secret communication in 6G networks \cite{Kihero2024,Chorti2022}, also in RIS-aided networks \cite{Khoshafa2024}. Some previous studies have analyzed the physical layer security of RIS-aided wireless networks. The work~\cite{naeem2023security} analyzes the security challenges of RIS-based wireless networks. Critical security threats to 6G networks due to unlawful use of RISs are highlighted, and open research issues and viable prospects of RISs in 6G were put forth. Similarly, the work~\cite{kaur2024survey} explores the potential of RISs for physical layer security of wireless networks. Several PLS concepts, performance metrics, and applicability in advanced wireless networks are discussed. In addition, the work provides a comprehensive and systematic classification of RIS-aided wireless network topologies, including various system and channel fading models. The study~\cite{xu2023reconfiguring} emphasizes the importance of communication security in the 6G ecosystem and highlights the role of RISs in achieving physical-layer security, focusing on secrecy rate and secrecy outage probability. Motivated by these considerations, this work studies the problem of designing a wireless network aided by a reconfigurable intelligent surface in order to maximize both the energy efficiency and the security of the communication, analyzing the trade-off between active and passive RISs. 

\subsection{Prior Works}
EE and physical layer security can be considered together through the notion of secrecy energy efficiency (SEE), defined as the amount of bits that can be reliably and securely transmitted per Joule of consumed energy \cite{ZapJSTSP16}. In the context of RIS-based communications, most previous works consider the maximization of the system secrecy sum-rate (SSR), neglecting the network performance in terms of SEE. In~\cite{Li2023}, for instance, the authors present a non-orthogonal multiple access (NOMA) network employing a simultaneous transmission and reception (STAR) RIS. Additionally, the system secrecy outage probability is characterized, and the results compare favorably with related studies. Similarly, a NOMA-based network aided by a STAR-RIS is considered in~\cite{Pei2023}. Closed-form approximations of the secrecy outage probability are derived and evaluated using simulations. In~\cite{Trigui2021}, the secrecy outage and average rate of an RIS-aided network are evaluated under the assumption that the RIS is capable of applying discrete phase shifts. Similarly, analytical expressions for the ergodic SSR of an RIS-empowered wireless network with multiple eavesdroppers are derived~\cite{Xu2021}. The work~\cite{Zhang2022} considers a STAR-aided NOMA-based network, maximizing the worst-case secrecy capacity of the system. Furthermore, the ergodic secrecy capacity of an RIS-aided wireless network is analyzed and approximated in closed form, considering flying eavesdroppers~\cite{Wei2022}. In ~\cite{Yang2020}, the secrecy outage probability of an RIS-aided network is analyzed. In \cite{Zhao2023}, the secrecy rate of an RIS-aided network powered by wireless power transfer is optimized. In~\cite{Hoang2023}, the SSR of an RIS-aided network with space-ground communications is optimized. In~\cite{Yizhi2023}, the SSR of a multi-user RIS-aided wireless network is analyzed. In a similar fashion, the data interleaving method is employed to achieve secret communications, leveraging an RIS, in~\cite{Cai2023}. The SEE of RIS-based networks is considered in fewer works. In~\cite{Yichi2023}, the optimization of the SEE of an RIS-aided network is conducted using deep reinforcement learning. The work~\cite{Yang2023} employs alternating maximization and sequential programming to maximize the SEE of a multi-user network. Similarly, the authors of \cite{Hao23} employ a blend of sequential programming and alternating optimization to optimize the minimum SEE of a multi-user network. However, these last works also focus on nearly-passive RISs without addressing the use of active RISs. In \cite{Soderi2024}, it is shown that RISs can be used to provide secrecy and user authentication in visible light communication systems. 

The study of active RISs is mostly confined to the reliability of wireless communications, without considering the secrecy aspect. Moreover, most available studies that consider energy aspects in RIS-aided networks focus on minimizing the network power consumption rather than considering the more general problem of EE maximization. In \cite{Chen2022}, an active omnidirectional RIS is used to improve the performance of vehicular networks. The base station (BS) precoding matrix and the RIS reflection coefficients are optimized to minimize the BS transmit power, assuming imperfect channel knowledge. In \cite{Zhu2022}, the authors develop a novel sub-connected architecture for active RISs, in which a single amplifier drives multiple adjacent reflecting elements. Sum-rate maximization and power consumption minimization are tackled by optimizing the RIS reflection coefficients, and BS transmit powers. The tradeoff between fully connected and sub-connected architectures for active RISs is investigated in \cite{Liu2022}, with reference to a multiple-antenna system. The maximization of the network EE is carried out with respect to the RIS reflection coefficients and transmit beamforming. In \cite{Ma2023}, the EE of a multi-user MISO system aided by an active RIS is optimized with respect to the BS beamforming and the RIS reflection coefficients. The problem is tackled by means of the quadratic transform method for fractional problems \cite{Shen2018}. A recent study in \cite{Peng2023} considers a hybrid RIS with both nearly-passive and active elements and aims at optimizing the number of passive and active elements deployed to maximize the minimum among the ergodic energy efficiencies of the mobile users, considering a single-antenna BS. In \cite{Fotock2023,FotockICC24}  the EE of a multi-user wireless network is considered, and resource allocation algorithms are provided to optimize the RIS reflection coefficients, the users transmit powers, and BS receive filters. In \cite{Lv2023}, the problem of minimizing the transmit power subject to a minimum secrecy rate constraint is tackled by optimizing the transmit beamforming and the RIS reflection matrix in a single-user wire-tap channel. In this context, it is shown that an active RIS can effectively combat the multiplicative fading effect while providing satisfactory security and power consumption levels. In \cite{Wang2022}, satellite communication in the presence of an eavesdropper is considered, in which an active RIS with local reflection capabilities is used to boost the received power. The minimum of the mobile users’ SEEs is maximized with respect to the RIS reflection coefficients and transmit beamforming by an alternating maximization technique, assuming a bounded error model for the channel between the RIS and the eavesdropper. 
The work~\cite{naeem2023security} analyzes the security challenges of RIS-based wireless networks. Critical security threats to 6G networks due to unlawful use of RISs are highlighted, and open research issues and viable prospects of RISs in 6G were broached. Similarly, the work~\cite{kaur2024survey} explores the potential of RISs for physical layer security of wireless networks. Several PLS concepts, performance metrics, and applicability in advanced wireless networks are discussed. In addition, the work provides a comprehensive and systematic classification of RIS-aided wireless network topologies, including various system and channel fading models. The study~\cite{xu2023reconfiguring} emphasizes the importance of communication security in the 6G ecosystem and highlights the role of RISs in achieving physical-layer security, focusing on secrecy rate and secrecy outage probability.

\subsection{Contributions}
In order to fill the aforementioned gaps, this work considers the uplink of a cellular network in which mobile users reach the BS with the aid of an RIS, while securing the content of the communication against the presence of a node that might eavesdrop the communication. The work focuses on the analysis and optimization of the system radio resources in order to maximize the SEE. The following contributions are made:
\begin{itemize}
	\item Two novel resource allocation algorithms are developed to optimize the users' transmit powers, the reflection coefficients at the RIS, and the base station receive filters for SEE maximization, thus addressing at the same time the energy and security aspects of the wireless network. 
	\item The optimization is carried out considering both scenarios of perfect channel state information (CSI) of all wireless channels and assuming only statistical knowledge of the channel to the eavesdropper. 
	\item The models used for the problem formulation assume an active RIS, but they encompass the use of nearly-passive RISs as a special case. This allows analyzing the trade-off between active and nearly-passive RISs in terms of SEE.
	\item Unlike most previous studies, this work considers the use of RISs with global reflection coefficients. These are a new type of RISs that generalizes traditional RISs with local reflection capabilities by enforcing the constraint on the reflected power considering the whole surface instead of each individual reflecting element \cite{Renzo2022}. 
\end{itemize}
The rest of this work is organized as follows. Section \ref{Sec:SysModel} defines the system model and problem formulation, considering perfect and statistical channel state information (CSI) for  RISs with global reflections. Section \ref{Sec:OPT_PCSI} presents the optimization with perfect CSI. Similarly, Section \ref{Sec:OPT_SCSI} covers the optimization with statistical CSI. Section \ref{Sec:NUM_ANA} provides the numerical analysis and discussion of results. Finally, conclusions are given in Section \ref{Sec:Concl}

\subsection{Notation}
$\mathbb{C}^{M \times N}$ and $\mathbb{R}^{M \times N}$ represent the spaces of $M \times N$ complex and real matrices, respectively. Scalars, column vectors, and matrices are denoted by lowercase, boldface lowercase, and boldface uppercase letters, e.g. $x$, $\mathbf{x}$, $\mathbf{X}$. $|x|$ represents the modulus of a scalar $x$, while $\|\mathbf{x}\|$ denotes the Euclidean norm of a column vector $\mathbf{x}$.  $(\cdot)^*$, $(\cdot)^T$, $(\cdot)^H$ indicate the conjugate, transpose, and Hermitian operations, respectively. $\text{diag}(\cdot)$ denotes a diagonal matrix.  $tr(\mathbf{S})$ and $\mathbf{S}^{-1}$ denote the trace and inverse of a square matrix. $\mathbf{S} \succeq \mathbf{0}$ indicates that $\mathbf{S}$ is positive semi-definite. $\mathbf{I}_M$ denotes an identity matrix of size $M \times M$. Statistical expectations are denoted by $\mathbb{E}(\cdot)$. For complex numbers, the real and imaginary parts are denoted as $\Re\{\cdot\}$ and $\Im\{\cdot\}$, respectively. The distribution of a circularly symmetric complex Gaussian (CSCG) random variable $x$ with a mean of $\mu$ and a variance of $\sigma^2$ is represented as $x \sim \mathcal{CN}(\mu, \sigma^2)$.
	
\section{System Model and Problem Formulation}\label{Sec:SysModel}
\subsection{System model}
\begin{figure}[!h]
	\centering
	\includegraphics[width=0.5\textwidth]{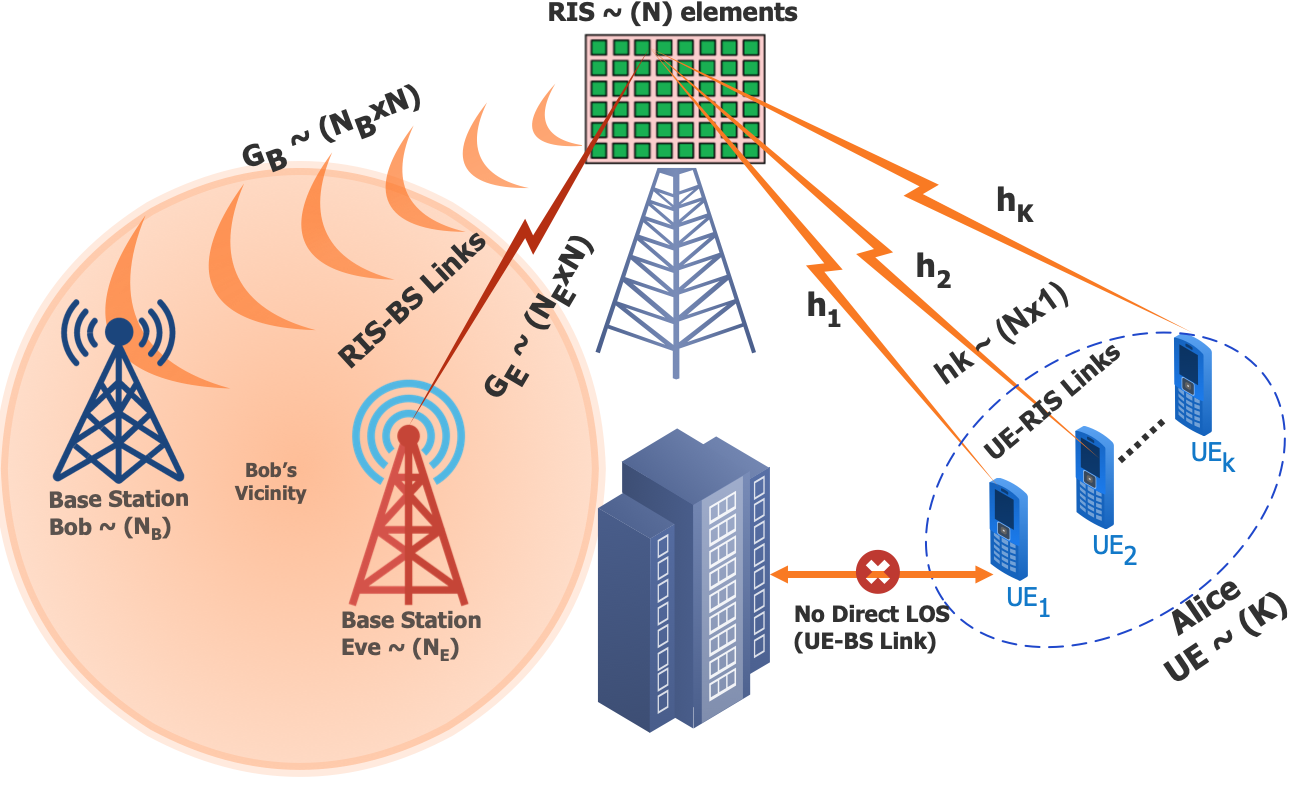}\caption{RIS-aided wireless network scenario} \label{Fig:scenario}
\end{figure}

We consider a network consisting of $K$ single-antenna transmitters, which communicate with their intended receiver, equipped with $N_{B}$ antennas, through an active RIS, equipped with $N$ reflecting elements. In the same area, we assume the presence of a single-antenna eavesdropper, which also receives the signal reflected by the RIS. The eavesdropper can be another mobile node of the network that is present in the area and inevitably overhears the communication or a malicious eavesdropper. The considered scenario is depicted in Fig.~\ref{Fig:scenario}. 

\noindent Let us denote by $p_{k}$ the $k$-th user's transmit power, by $\boldsymbol{\gamma} = \left(\gamma_{1},\ldots,\gamma_{N}\right)$ the $N \times 1$ vector containing the RIS reflection coefficients, by $\boldsymbol{h}_{k}$ the $N \times 1 $ channel between the $k$-th user and the RIS, by $\boldsymbol{G}_{B}$ the $N_{B} \times N$ channel from the RIS to the legitimate receiver, and by $\boldsymbol{g}$ the $N \times 1$ channel from the RIS to the eavesdropper. In this setup, the SINR enjoyed by user $k$ at the intended and eavesdropping receivers can be written as 
\begin{align}\label{Eq:SINR_B}
	\text{SINR}_{k,B}&=\frac{p_{k}\left|\boldsymbol{c}_{k}^{H}\boldsymbol{A}_{k}\boldsymbol{\gamma}\right|^{2}}{\boldsymbol{c}_{k}^{H}\boldsymbol{W}\boldsymbol{c}_{k}+\sum_{m\neq k}p_{m}\left|\boldsymbol{c}_{k}^{H}\boldsymbol{A}_{m}\boldsymbol{\gamma}\right|^{2}}\\
\label{Eq:SINR_E}
	\text{SINR}_{k,E}&=\frac{p_{k}\left|\boldsymbol{g}^{H}\boldsymbol{H}_{k}\boldsymbol{\gamma}\right|^{2}}{\sigma^{2}_{E}+\sigma_{RIS}^{2}\boldsymbol{g}^{H}\boldsymbol{\Gamma}\boldsymbol{\Gamma}^{H}\boldsymbol{g}+\sum_{m\neq k}p_{m}\left|\boldsymbol{g}^{H}\boldsymbol{H}_{m}\boldsymbol{\gamma}\right|^{2}}\;,
\end{align}
where $\boldsymbol{A}_{k}=\boldsymbol{G}_{B}\boldsymbol{H}_{k}$, with $\boldsymbol{H}_{k}=\text{diag}(\boldsymbol{h}_{k})$, while $\boldsymbol{c}_{k}$ is the linear receive filter applied by the intended receiver to decode the data from user $k$, $\boldsymbol{W} =\sigma_{B}^{2}\boldsymbol{I}_{N_{B}}+ \sigma_{RIS}^{2}\boldsymbol{G}_{B}\boldsymbol{\Gamma}\boldsymbol{\Gamma}^{H}\boldsymbol{G}_{B}^{H}$ is the overall covariance matrix at the intended receiver, wherein $\boldsymbol{\Gamma} = \text{diag}\left(\boldsymbol{\gamma}\right)$ and $\sigma^{2}_{RIS}$ the noise variance introduced at the active RIS. In this context, the system SSR is expressed as \footnote{We assume that the channel and propagation condition are such that feasible system operating points exist at which each user experiences a non-negative secrecy rate. Otherwise, the scenario would not be of practical interest from a physical layer security perspective.}
\begin{align}\label{Eq:SecRate}
	\text{SSR}&=\sum_{k=1}^{K}\log_{2}\left(1+\text{SINR}_{k,B}\right) - \log_{2}\left(1+\text{SINR}_{k,E}\right).
\end{align}
In the following, we develop a model for both active and nearly-passive RISs. The difference between the two is that the former are equipped with an analog amplifier which can increase the radio-frequency power of the input signal, while the latter do not have any amplification capability and require only the static power that is needed to operate the reconfiguration circuits. Moreover, we consider RISs with global reflection capabilities. These generalize RISs with local reflection capabilities, in which the reflection constraints are applied separately to each reflector, i.e. $|\gamma(n)|^{2}\leq 1$ for all $n=1,\ldots,N$ for nearly-passive RISs and $|\gamma(n)|^{2}\leq P_{r,max}$ for all $n=1,\ldots,N$ for active RISs with $P_{r,max}$ the maximum radio-frequency power that can be provided by the power amplifier. In an RIS with global reflection capabilities, a single reflection constraint applies to all of the coefficients $\{\gamma(n)\}_{n=1}^{N}$. If the RIS is nearly-passive, the constraint must ensure that $P_{out}\leq P_{in}$, i.e. the total radio-frequency power reflected by the whole RIS must be not larger than the total radio-frequency power that impinges on the RIS. If instead the RIS is active, the global reflection constraint is $P_{out}\leq P_{in}+P_{r,max}$, which ensures that the output power must not be larger than the total input power plus the maximum radio-frequency power that can be provided by the RIS amplifier $P_{r,max}$. In the following, we will consider explicitly the case of active RIS with global reflection constraints, which subsumes the nearly-passive scenario as a special case, as will be discussed in Remark \ref{Rem:PowerCons}. 

To elaborate, the radio-frequency power consumed by the active RIS is given by the difference between the output and input power, which can be computed as \cite{fotock2023energy}
\begin{align}\label{Eqn:P_RF}
	P_{out} - P_{in} &= \tr\left(\sum_{k=1}^{K}p_{k}\boldsymbol{\Gamma}\boldsymbol{h}_{k}\boldsymbol{h}^{H}_{k}\boldsymbol{\Gamma}^{H} +  \sigma^{2}_{RIS}\boldsymbol{\Gamma}\boldsymbol{\Gamma}^{H}\right)\\
	&- \sum_{k=1}^{K}p_{k} \left\Vert \boldsymbol{h}_{k}\right\Vert^{2} - \sigma^{2}_{RIS}N = \tr\left(\left(\boldsymbol{\gamma}\boldsymbol{\gamma}^{H}-\boldsymbol{I}_{N}\right)\boldsymbol{R}\right)\notag
\end{align}
wherein $\boldsymbol{R}=\textstyle\sum_{k=1}^{K}p_{k}\boldsymbol{H}_{k}^{H}\boldsymbol{H}_{k}+\sigma_{RIS}^{2}\boldsymbol{I}_{N}$. Thus, the total power consumption of the system is given by 
\begin{equation}\label{Eq:Power}
P_{tot}=\tr\left(\left(\boldsymbol{\gamma}\boldsymbol{\gamma}^{H}-\boldsymbol{I}_{N}\right)\boldsymbol{R}\right)+\sum_{k=1}^{K}p_{k}+P_{c}\;,
\end{equation}
wherein $P_{c} = NP_{c,n} + P_{0,RIS} + P_{0}$, with $P_{c,n} $ the static power consumption of the $n$-th RIS element, $P_{0,RIS}$ is the rest of the static power consumed by the RIS and $P_{0}$ accounts for the static power consumption of the transmitters and legitimate receiver. 
\begin{remark}\label{Rem:PowerCons}
The power consumption in \eqref{Eq:Power} has been developed for an active RIS, but it specializes to the case of a nearly-passive RIS by simply removing the term $P_{out} - P_{in}$. Indeed, for nearly-passive RISs, $P_{out} \leq P_{in}$ and, thus, the difference $P_{out} - P_{in}$ does not appear in the power consumption $P_{tot}$. Moreover, the terms $P_{c,n}$ and $P_{0,RIS}$ will have a lower numerical value than in the case of an active RIS because nearly-passive RISs require simpler hardware.
\end{remark}
Thus, the SEE is given by the ratio between the SSR and the network power consumption, namely
\begin{equation}\label{Eq:SEE_PCSI}
\text{SEE}=\frac{\log_{2}\left(1+\text{SINR}_{B}\right) - \log_{2}\left(1+\text{SINR}_{E}\right)}{\tr\left(\left(\boldsymbol{\gamma}\boldsymbol{\gamma}^{H}-\boldsymbol{I}_{N}\right)\boldsymbol{R}\right)+\sum_{k=1}^{K}\mu_{k}p_{k}+P_{c}}\;,
\end{equation}

\subsection{Problem formulation}
The rest of this work is concerned with maximizing the SEE, which will be performed considering two CSI assumptions. 
\subsubsection{Perfect CSI} we assume that all wireless channels are estimated with negligible error at the legitimate receiver. As for the channels between the RIS and the transmitters/legitimate receivers, perfect CSI can be obtained by standard channel estimation methods. As for the channel between the RIS and the eavesdropper, the assumption of perfect CSI applies to a scenario in which the eavesdropper is not a malicious, hidden node but rather a terminal whose presence in the coverage area is allowed, e.g., a node at the cell-edge that is served by a neighboring BS. In this situation, the eavesdropper inevitably receives the signal intended for the legitimate receiver and thus could decode the information sent by the transmitters, which should be avoided for privacy reasons. In this context, the SEE maximization problem can be formulated as
\begin{subequations}\label{Prob:SEE_ARIS}
	\begin{align}
		&\ds\max_{\boldsymbol{\gamma},\boldsymbol{p},\boldsymbol{C}}\; \text{SEE}(\boldsymbol{\gamma},\boldsymbol{p},\boldsymbol{C})\label{Prob:aSEE_ARIS}\\
		&\;\text{s.t.}\;\tr\left(\boldsymbol{R}\right)\leq \tr(\boldsymbol{R}\boldsymbol{\gamma}\boldsymbol{\gamma}^{H})\leq P_{R,max}+\tr\left(\boldsymbol{R}\right)\label{Prob:bSEE_ARIS}\\
		&\;\quad\;\;0\leq p_{k}\leq P_{max,k}\;\forall\;k=1,\ldots,K\;,\label{Prob:cSEE_ARIS}
	\end{align}
\end{subequations}
with $\boldsymbol{C} = \left[\boldsymbol{c}_{1},\ldots,\boldsymbol{c}_{K} \right]$.\footnote{We observe that Problem~\eqref{Prob:SEE_ARIS} is always feasible, since setting $\left|\gamma_{n}\right|=1$ for all $n$ fulfills all constraints.}

\subsubsection{Statistical CSI} we assume that all wireless channels, except the channel between the RIS and the eavesdropper, are estimated with negligible error at the legitimate receiver. Instead, as for the channel between the RIS and the eavesdropper, we assume that only statistical CSI is available. Specifically, we follow a mean feedback model, assuming that the true channel is given by $\boldsymbol{g}=\boldsymbol{\widehat{g}}+\boldsymbol{\delta}$, wherein $\boldsymbol{\widehat{g}}$ is the mean value of $\boldsymbol{g}$, which is estimated by the legitimate receiver, but subject to the random error $\boldsymbol{\delta}$, for which it is only known that $\boldsymbol{\delta}\sim {\cal CN}(\boldsymbol{0},\sigma_{g}^{2}\boldsymbol{I}_{N})$. This scenario applies to a situation in which the eavesdropper is a hidden node whose position is known up to some estimation error. As an example, this scenario applies to an indoor setting in which a possible eavesdropper might be located in a given room or floor. Moreover, the parameter $\sigma_{g}^{2}$ controls the reliability of the channel estimate $\boldsymbol{\widehat{g}}$. The limiting case $\sigma_{g}^{2}>>\|\boldsymbol{\widehat{g}}\|^{2}$, models the scenario in which $\boldsymbol{g}$ is completely unknown, while we fall back to the perfect CSI case by setting $\sigma_{g}^{2}=0$. If only statistical CSI case is available, the SEE maximization problem can be formulated as
\begin{subequations}\label{Prob:SEE_ARIS_SCSI}
	\begin{align}
		&\ds\max_{\boldsymbol{\gamma},\boldsymbol{p},\boldsymbol{C}}\; \mathbb{E}_{\boldsymbol{\delta}}\left[\text{SEE}(\boldsymbol{\gamma},\boldsymbol{p},\boldsymbol{C})\right]\label{Prob:aSEE_ARIS_SCSI}\\
		&\;\text{s.t.}\;\tr\left(\boldsymbol{R}\right)\leq \tr(\boldsymbol{R}\boldsymbol{\gamma}\boldsymbol{\gamma}^{H})\leq P_{R,max}+\tr\left(\boldsymbol{R}\right)\label{Prob:bSEE_ARIS_SCSI}\\
		&\;\quad\;\;0\leq p_{k}\leq P_{max,k}\;\forall\;k=1,\ldots,K\;,\label{Prob:cSEE_ARIS_SCSI}
	\end{align}
\end{subequations}
wherein we have taken the average of the SEE with respect to the random vector $\boldsymbol{\delta}$.\footnote{Note that the constraints are not affected by the 
randomness of the channel between the RIS and the eavesdropper.}

\begin{remark}\label{Rem:Extension}
The model developed in this section can be readily extended to consider multi-cell networks, in which interference is treated as noise. Indeed, in this case it suffices to include the out-of-cell interference in the sum over the interfering user, which formally does not change the mathematical expression of either $\text{SINR}_{B}$, $\text{SINR}_{E}$, or $P_{tot}$. Therefore, the algorithms to be developed in the rest of this document can also be generalized to work in multi-cell setups. 
\end{remark}

\section{Optimization with Perfect CSI}\label{Sec:OPT_PCSI}
This section tackles Problem \eqref{Prob:SEE_ARIS} by the popular alternating maximization method, optimizing cyclically the RIS reflection vector $\boldsymbol{\gamma}$, the users' transmit powers $\boldsymbol{p}$, and the linear receive filters $\boldsymbol{C}$. These three subproblems will be treated separately. 

\subsection[Optimization of gamma]{Optimization of $\boldsymbol{\gamma}$}
The optimization with respect to the RIS vector $\boldsymbol{\gamma}$ is 
\begin{subequations}\label{Prob:SEC_ARIS_gamma_PCSI}
	\begin{align}
		&\ds\max_{\boldsymbol{\gamma}}\; B\frac{\textstyle\sum_{k=1}^{K}\log_{2}\left(1+\text{SINR}_{k,B}\right) - \log_{2}\left(1+\text{SINR}_{k,E}\right)}{\tr\left(\boldsymbol{R}\boldsymbol{\gamma}\boldsymbol{\gamma}^{H}\right)+P_{c,eq}} \label{Prob:aSEC_ARIS_gamma_PCSI}\\
		&\;\text{s.t.}\;\tr\left(\boldsymbol{R}\right)\leq \tr(\boldsymbol{R}\boldsymbol{\gamma}\boldsymbol{\gamma}^{H})\leq P_{R,max}+\tr\left(\boldsymbol{R}\right)\label{Prob:bSEC_ARIS_gamma_PCSI}
	\end{align}
\end{subequations}
wherein $P_{c,eq}=\sum_{k}\mu_{k}p_{k}+P_{c}-\tr(\boldsymbol{R})$. The objective in \eqref{Prob:aSEC_ARIS_gamma_PCSI} is neither concave nor pseudo-concave in $\boldsymbol{\gamma}$. Moreover, the first inequality in~\eqref{Prob:bSEC_ARIS_gamma_PCSI} is a non-convex constraint. Thus, directly applying fractional programming is not feasible~\cite{ZapNow15}. To address this challenge, we employ the sequential fractional programming method. As a first step, we express the term $\boldsymbol{c}_{k}^{H}\boldsymbol{W}\boldsymbol{c}_{k}$ in terms of the vector $\boldsymbol{\gamma}$, instead of the matrix $\boldsymbol{\Gamma}$. To this end, we define $\boldsymbol{u}_{k}=\boldsymbol{G}^{H}\boldsymbol{c}_{k}$ and $\widetilde{\boldsymbol{U}}_{k}=\text{diag}\left(|u_{1}|^{2},\ldots,|u_{N}|^{2}\right)$. Hence, we obtain 
\begin{align}\label{Eq:W_gamma}
\boldsymbol{c}_{k}^{H}\boldsymbol{W}\boldsymbol{c}_{k}&=
\sigma^{2}\|\boldsymbol{c}_{k}\|^{2}+\sigma_{RIS}^{2}\boldsymbol{\gamma}^{H}\widetilde{\boldsymbol{U}}_{k}\boldsymbol{\gamma}\;.
\end{align}
Moreover, we also observe that $\boldsymbol{g}^{H}\boldsymbol{\Gamma}\boldsymbol{\Gamma}^{H}\boldsymbol{g}=\|\boldsymbol{G}\boldsymbol{\gamma}\|^{2}$, wherein $\boldsymbol{G}=\text{diag}(g_{1},\ldots,g_{N})$, with $g_{n}$ the $n$-th component of the vector $\boldsymbol{g}$. Next, defining
\begin{align}\label{Eq:X_B}
	&x_{B}=p_{k}\left|\boldsymbol{c}_{k}^{H}\boldsymbol{A}_{k}\boldsymbol{\gamma}\right|^{2}\;,\; x_{E}=p_{k}|\boldsymbol{g}^{H}\boldsymbol{H}_{k}\boldsymbol{\gamma}|^{2}\\
	& y_{B}=\sigma^{2}\|\boldsymbol{c}_{k}\|^{2}+\sigma_{RIS}^{2}\|\widetilde{\boldsymbol{U}}_{k}^{1/2}\boldsymbol{\gamma}\|^{2}+\textstyle\sum_{m\neq k}p_{m}\left|\boldsymbol{c}_{k}^{H}\boldsymbol{A}_{m}\boldsymbol{\gamma}\right|^{2}\label{Eq:Y_B}\\
	&y_{E}=\textstyle\sum_{m\neq k}p_{m}|\boldsymbol{g}^{H}\boldsymbol{H}_{m}\boldsymbol{\gamma}|^{2}+\sigma_{RIS}^{2}\|\boldsymbol{G}\boldsymbol{\gamma}\|^{2}\;.\label{Eq:Y_E}
\end{align}
The $k$-th summand of the numerator of \eqref{Prob:aSEC_ARIS_gamma_PCSI} can be written as 
\begin{align}\label{Eq:SecRate_PCSI}
	R_{s,k}^{PCSI}&\!=\!\log_{2}\!\left(\!1\!+\!\frac{x_{B}}{y_{B}}\!\right)\!-\!\log_{2}\!\left(\!1\!+\!\frac{x_{E}}{y_{E}\!+\!\sigma_{E}^{2}}\right)\!=\!\log_{2}\!\left(\!1\!+\!\frac{x_{B}}{y_{B}}\!\right)\notag\\
	&\!+\!\log_{2}\left(1\!+\!\frac{y_{E}}{\sigma_{E}^{2}}\right)\!+\!\log_{2}\left(\sigma_{E}^{2}\right)\!-\!\log_{2}(\sigma_{E}^{2}\!+\!x_{E}\!+\!y_{E}).
\end{align}
In order to lower-bound \eqref{Eq:SecRate_PCSI}, we will leverage the following lower-bounds \begin{align}\label{Eq:LB1}
	\log_{2}\left(1\!+\!\frac{x}{y}\right)\!\geq\! \log_{2}\left(1\!+\!\frac{\bar{x}}{\bar{y}}\right)\!+\!\frac{\bar{x}}{\bar{y}}\left(\frac{2\sqrt{x}}{\sqrt{\bar{x}}}\!-\!\frac{x+y}{\bar{x}\!+\!\bar{y}}\!-\!1\right)\\
 \label{Eq:LB2}
	\log_{2}(\sigma_{E}^{2}\!+\!x\!+\!y)\leq \log_{2}(\sigma_{E}^{2}+\bar{x}\!+\!\bar{y})\!+\!\frac{1}{\ln{2}}\frac{x+y-\bar{x}-\bar{y}}{\sigma_{E}^{2}+\bar{x}+\bar{y}}\;,
\end{align}
which hold for any $x$, $y$, $\bar{x}$ and $\bar{y}$, and hold with equality whenever $x=\bar{x}$ and $y=\bar{y}$. Specifically, the first inequality lower-bounds the first and second term in \eqref{Eq:SecRate_PCSI}, while the second inequality lower-bounds the last term, accounting for the negative sign. To elaborate, for any feasible vector $\bar{\boldsymbol{\gamma}}$, let us define 
\begin{align}
	\hspace{-0.3cm}\bar{x}_{B}&=p_{k}\left|\boldsymbol{c}_{k}^{H}\boldsymbol{A}_{k}\boldsymbol{\bar{\gamma}}\right|^{2}\;,\;\bar{x}_{E}=p_{k}|\boldsymbol{g}^{H}\boldsymbol{H}_{k}\boldsymbol{\bar{\gamma}}|^{2}\label{Eq:X_B_bar}\\
	\hspace{-0.3cm}\bar{y}_{B}&=\sigma^{2}\|\boldsymbol{c}_{k}\|^{2}+\sigma_{RIS}^{2}\|\widetilde{\boldsymbol{U}}_{k}^{1/2}\boldsymbol{\bar{\gamma}}\|^{2}+\textstyle\sum_{m\neq k}p_{m}\left|\boldsymbol{c}_{k}^{H}\boldsymbol{A}_{m}\boldsymbol{\bar{\gamma}}\right|^{2}\label{Eq:Y_B_bar}\\
	\hspace{-0.3cm}\bar{y}_{E}&=\textstyle\sum_{m\neq k}p_{m}|\boldsymbol{g}^{H}\boldsymbol{H}_{m}\boldsymbol{\bar{\gamma}}|^{2}+\sigma_{RIS}^{2}\|\boldsymbol{G}\boldsymbol{\bar{\gamma}}\|^{2}\;,\label{Eq:Y_E_bar}
\end{align}
which yields the lower-bound
\begin{align}\label{Eq:LowerBound_PCSI}
	&R_{s,k}^{PCSI}(\boldsymbol{\gamma})\geq \log_{2}\left(1+\frac{\bar{x}_{B}}{\bar{y}_{B}}\right)\!+\!\frac{\bar{x}_{B}}{\bar{y}_{B}}\left(\frac{2\sqrt{x_{B}}}{\sqrt{\bar{x}_{B}}}-\frac{x_{B}+y_{B}}{\bar{x}_{B}+\bar{y}_{B}}-1\right)\notag\\
	&+\log_{2}\left(1+\frac{\bar{y}_{E}}{\sigma_{E}^{2}}\right)\!+\!\frac{\bar{y}_{E}}{\sigma_{E}^{2}}\left(\frac{2\sqrt{y_{E}}}{\sqrt{\bar{y}_E}}-\frac{y_{E}+\sigma_{E}^{2}}{\bar{y}_{E}+\sigma_{E}^{2}}-1\right)\notag\\
	&+\log_{2}\left(\frac{\sigma_{E}^{2}}{\sigma_{E}^{2}+\bar{x}_{E}+\bar{y}_{E}}\right)\!-\!\frac{x_{E}+y_{E}-\bar{x}_{E}-\bar{y}_{E}}{\sigma_{E}^{2}+\bar{x}_{E}+\bar{y}_{E}}=\widebar{R}_{s,k}^{PCSI}(\boldsymbol{\gamma}).
\end{align}
All terms in \eqref{Eq:LowerBound_PCSI} that depend on $\boldsymbol{\gamma}$ are concave, except for $\sqrt{x_{B}}$ and $\sqrt{y_{E}}$, which are convex. Indeed, $\sqrt{x_{B}}=\sqrt{p_{k}}\left|\boldsymbol{c}_{k}^{H}\boldsymbol{A}_{k}\boldsymbol{\gamma}\right|$, which is clearly convex in $\boldsymbol{\gamma}$, while $\sqrt{y_{E}}$ can be written as 
\begin{align}
\sqrt{y_{E}}&=\sqrt{\boldsymbol{\gamma}^{H}\left(\sum_{m\neq k}p_{m}\boldsymbol{H}_{m}^{H}\boldsymbol{g}\boldsymbol{g}^{H}\boldsymbol{H}_{m}+\sigma_{RIS}^{2}\boldsymbol{G}\boldsymbol{G}^{H}\right)\boldsymbol{\gamma}}\notag\\
&=\left\|\left(\sum_{m\neq k}p_{m}\boldsymbol{H}_{m}^{H}\boldsymbol{g}\boldsymbol{g}^{H}\boldsymbol{H}_{m}+\sigma_{RIS}^{2}\boldsymbol{G}\boldsymbol{G}^{H}\right)^{1/2}\boldsymbol{\gamma}\right\|\;,
\end{align}
which is convex, too. Thus, $\sqrt{x_{B}}$ and $\sqrt{y_{E}}$ can be lower-bounded by their first-order Taylor expansion around $\boldsymbol{\bar{\gamma}}$. Namely, defining the matrices
\begin{align}
\boldsymbol{M}_{k,B}&=\boldsymbol{A}_{k}^{H}\boldsymbol{c}_{k}\boldsymbol{c}_{k}^{H}\boldsymbol{A}_{k}\;,\\
\boldsymbol{M}_{k,E}&=\sum_{m\neq k}p_{m}\boldsymbol{H}_{m}^{H}\boldsymbol{g}\boldsymbol{g}^{H}\boldsymbol{H}_{m}+\sigma_{RIS}^{2}\boldsymbol{G}\boldsymbol{G}^{H}\;,
\end{align}
$\sqrt{x_{B}}$ and $\sqrt{y_{E}}$ can be lower-bounded by linear terms as: 
\begin{align}\label{Eq:lowerXB}
&\sqrt{x_{B}}=\sqrt{p_{k}}\sqrt{|\boldsymbol{c}_{k}^{H}\boldsymbol{A}_{k}\boldsymbol{\gamma}|^{2}}=\sqrt{p_{k}}\sqrt{\boldsymbol{\gamma}^{H}\boldsymbol{M}_{k,B}\boldsymbol{\gamma}}\geq\\
&\sqrt{p_{k}}\left(\sqrt{\boldsymbol{\bar{\gamma}}^{H}\boldsymbol{M}_{k,B}\boldsymbol{\bar{\gamma}}}+\frac{\Re{(\boldsymbol{\bar{\gamma}}^{H}\boldsymbol{M}_{k,B}}(\boldsymbol{\gamma}-\boldsymbol{\bar{\gamma}}))}{\sqrt{|\boldsymbol{c}_{k}^{H}\boldsymbol{A}_{k}\boldsymbol{\bar{\gamma}}|}}\right)=\widetilde{x}_{B}\;,\notag\\
&\sqrt{y_{E}}=\sqrt{\|\boldsymbol{M}_{k,E}^{1/2}\boldsymbol{\gamma}\|^{2}}=\sqrt{\boldsymbol{\gamma}^{H}\boldsymbol{M}_{k,E}\boldsymbol{\gamma}}\geq\label{Eq:lowerYE}\\
&\sqrt{\boldsymbol{\bar{\gamma}}^{H}\boldsymbol{M}_{k,E}\boldsymbol{\bar{\gamma}}}+\frac{\Re{(\boldsymbol{\bar{\gamma}}^{H}\boldsymbol{M}_{k,E}}(\boldsymbol{\gamma}-\boldsymbol{\bar{\gamma}}))}{\sqrt{|\boldsymbol{c}_{k}^{H}\boldsymbol{A}_{k}\boldsymbol{\bar{\gamma}}|}}=\widetilde{y}_{E}\;. \notag
\end{align}
Then, plugging \eqref{Eq:lowerXB} and \eqref{Eq:lowerYE} into \eqref{Eq:LowerBound_PCSI} yields
\begin{align}
&R_{s,k}^{PCSI}(\boldsymbol{\gamma})\geq \widebar{R}_{s,k}^{PCSI}(\boldsymbol{\gamma}) \geq 
\log_{2}\left(1+\frac{\bar{x}_{B}}{\bar{y}_{B}}\right)\notag\\
&+\frac{\bar{x}_{B}}{\bar{y}_{B}}\left(\frac{2\widetilde{x}_{B}}{\sqrt{\bar{x}_{B}}}-\frac{x_{B}+y_{B}}{\bar{x}_{B}+\bar{y}_{B}}-1\right)+\log_{2}\left(1+\frac{\bar{y}_{E}}{\sigma_{E}^{2}}\right)\notag\\
&+\frac{\bar{y}_{E}}{\sigma_{E}^{2}}\left(\frac{2\widetilde{y}_{E}}{\sqrt{\bar{y}_E}}-\frac{y_{E}+\sigma_{E}^{2}}{\bar{y}_{E}+\sigma_{E}^{2}}-1\right)+\log_{2}\left(\frac{\sigma_{E}^{2}}{\sigma_{E}^{2}+\bar{x}_{E}+\bar{y}_{E}}\right)\notag\\
&-\frac{x_{E}\!+\!y_{E}\!-\!\bar{x}_{E}\!-\!\bar{y}_{E}}{\sigma_{E}^{2}\!+\!\bar{x}_{E}\!+\!\bar{y}_{E}}\!=\!\widetilde{R}_{s,k}^{PCSI}(\boldsymbol{\gamma})\;.
\end{align}
Lastly, it remains to deal with the non-convexity of the left-hand side of \eqref{Prob:bSEC_ARIS_gamma_PCSI}. This can be accomplished by resorting again to the sequential approximation framework. Specifically, since $\tr\left(\boldsymbol{R}\boldsymbol{\gamma}\boldsymbol{\gamma}^{H}\right)$ is convex in $\boldsymbol{\gamma}$, its first-order Taylor expansion around any point $\bar{\boldsymbol{\gamma}}$ provides a lower-bound:
\begin{align}\label{Eq:ConstraintSFP}
	\tr\left(\boldsymbol{R}\boldsymbol{\gamma}\boldsymbol{\gamma}^{H}\right) &= \boldsymbol{\gamma}\boldsymbol{R}\boldsymbol{\gamma}^{H} \geq \bar{\boldsymbol{\gamma}}\boldsymbol{R}\bar{\boldsymbol{\gamma}}^{H} + 2\Re\left\{\bar{\boldsymbol{\gamma}}^{H}\boldsymbol{R}\left(\boldsymbol{\gamma} - \bar{\boldsymbol{\gamma}}\right)\right\}
\end{align}
Consequently, in each iteration of the sequential method, the following problem must be solved:
\begin{subequations}\label{Prob:SEE_Gamma_Active_PCSI}
	\begin{align}
		&\ds\max_{\boldsymbol{\gamma}}\;\frac{\textstyle\sum_{k=1}^{K}\widetilde{R}_{s,k}}{\tr\left(\boldsymbol{R}\boldsymbol{\gamma}\boldsymbol{\gamma}^{H}\right)+P^{\text{(a)}}_{c,eq}} \label{Prob:aSEE_Gamma_Active_PCSI}\\
		&\;\text{s.t.}\;\boldsymbol{\gamma}\boldsymbol{R}\boldsymbol{\gamma}^{H} \leq P_{R,max}+\tr\left(\boldsymbol{R}\right)\label{Prob:bSEE_Gamma_Active_PCSI}\\
		&\; \quad \;\bar{\boldsymbol{\gamma}}\boldsymbol{R}\bar{\boldsymbol{\gamma}}^{H} + 2\Re\left\{\bar{\boldsymbol{\gamma}}^{H}\boldsymbol{R}\left(\boldsymbol{\gamma} - \bar{\boldsymbol{\gamma}}\right)\right\} \geq \tr\left(\boldsymbol{R}\right)\label{Prob:cSEE_Gamma_Active_PCSI}
	\end{align}
\end{subequations}
The objective in \eqref{Prob:SEE_Gamma_Active_PCSI} has a concave numerator and a convex denominator. Thus, \eqref{Prob:SEE_Gamma_Active_PCSI} is a pseudo-concave maximization with convex constraints, which can be solved by fractional programming. The procedure for RIS optimization is stated in Algorithm~\ref{Alg:SCA_gamma1_RIS}.
\begin{algorithm}
	\caption{RIS optimization - Perfect CSI}
	\begin{algorithmic}\label{Alg:SCA_gamma1_RIS}
		\STATE $\epsilon > 0$, $\boldsymbol{\bar{\gamma}}$ \texttt{any feasible vector};
		\REPEAT
		\STATE \texttt{Solve} \eqref{Prob:SEE_Gamma_Active_PCSI} \texttt{and let} $\boldsymbol{\gamma^{\star}}$ \texttt{be the solution};
        \STATE $T=|\text{SEE}(\boldsymbol{\gamma^{\star}})-\text{SEE}(\boldsymbol{\bar{\gamma}})|$; 
        \STATE $\bar{\boldsymbol{\gamma}}=\boldsymbol\gamma^{\star}$; 
		\UNTIL{$T<\epsilon$}
	\end{algorithmic}
\end{algorithm}

\begin{proposition}\label{Prop:SCA_gamma1_RIS}
	Algorithm \ref{Alg:SCA_gamma1_RIS} monotonically improves the value of the objective of Problem \eqref{Prob:SEC_ARIS_gamma_PCSI} and converges to a point fulfilling the Karush-Kuhn-Tucker (KKT) optimality conditions of \eqref{Prob:SEC_ARIS_gamma_PCSI}.
\end{proposition}
\begin{IEEEproof}
	The proof follows upon noticing that Algorithm~\ref{Alg:SCA_gamma1_RIS} fulfills all assumptions of the sequential fractional programming framework \cite{SeqCvxProg78}. Indeed, \eqref{Prob:aSEE_Gamma_Active_PCSI} is a lower-bound of the original objective in \eqref{Prob:aSEC_ARIS_gamma_PCSI}, which is tight at $\boldsymbol{\bar{\gamma}}$. Similarly, the left-hand-side of  \eqref{Prob:cSEE_Gamma_Active_PCSI} is a lower-bound of the left-hand-side in the first inequality of \eqref{Prob:bSEC_ARIS_gamma_PCSI}, and the bound is tight at $\boldsymbol{\bar{\gamma}}$. Finally, it can be easily verified that the first-order derivative of \eqref{Prob:aSEE_Gamma_Active_PCSI} and of the left-hand-side of \eqref{Prob:cSEE_Gamma_Active_PCSI} are equal, when evaluated at $\boldsymbol{\bar{\gamma}}$, to the first-order derivatives of \eqref{Prob:aSEC_ARIS_gamma_PCSI} and of the left-hand-side in the first inequality of \eqref{Prob:bSEC_ARIS_gamma_PCSI}, respectively. 
\end{IEEEproof}


\subsection[Optimization of transmit powers]{Optimization of $\boldsymbol{p}$}\label{Sec:PowerPCSI}
Let us define $a^{(B)}_{k,m} = |\boldsymbol{c}_{k}^{H}\boldsymbol{A}_{m}\boldsymbol{\gamma}|^{2}$ for all $m$ and $k$, $d_{k}^{(B)} = \boldsymbol{c}_{k}^{H}\boldsymbol{W}\boldsymbol{c}_{k}$, $a^{(E)}_{k} = |\boldsymbol{g}^{H}\boldsymbol{H}_{k}\boldsymbol{\gamma}|^{2}$, $d^{(E)} = \sigma_{E}^{2}+\sigma_{RIS}^{2}\boldsymbol{g}^{H}\boldsymbol{\Gamma}\boldsymbol{\Gamma}^{H}\boldsymbol{g}$, $\mu_{k}=1+\tr{\left((\boldsymbol{\gamma}\boldsymbol{\gamma}^{H}-\boldsymbol{I}_{N})\boldsymbol{H}_{k}^{H}\boldsymbol{H}_{k}\right)}$, and $P_{c,eq} = \sigma_{RIS}^{2}\left(\|\boldsymbol{\gamma}\|^{2}-N\right) + P_{c}$. Thus, the power optimization problem can be formulated as:
\begin{subequations}\label{Prob:SEE_ARIS_power_PCSI}
	\begin{align}
		&\ds\max_{\boldsymbol{p}}\,\frac{\sum_{k=1}^{K}\ds\log_{2}\left(1+\frac{p_{k}a^{(B)}_{k,k}}{d_{k}^{(B)}+\sum_{m\neq k}p_{m}a^{(B)}_{k,m}}\right)}{\sum_{k=1}^{K}\mu_{k}p_{k}+P_{c,eq}}\notag\\
		&-\frac{\sum_{k=1}^{K}\ds\log_{2}\left(1+\frac{p_{k}a^{(E)}_{k}}{d^{(E)}+\sum_{m\neq k}p_{m}a^{(E)}_{m}}\right)}{\sum_{k=1}^{K}\mu_{k}p_{k}+P_{c,eq}}\label{Prob:aSEE_ARIS_power_PCSI}\\
		&\;\text{s.t.}\;0\leq p_{k}\leq P_{max,k}\;,\forall\; k=1,\ldots,K.
	\end{align}
\end{subequations}
The objective \eqref{Prob:aSEE_ARIS_power_PCSI} is neither concave nor pseudo-concave due to the non-concavity of its numerator. As a result, fractional programming approaches can not solve Problem \eqref{Prob:aSEE_ARIS_power_PCSI} with affordable complexity. For this reason, we employ the sequential fractional programming method \cite{ZapNow15}, which can find a first-order optimal solution of \eqref{Prob:aSEE_ARIS_power_PCSI} with polynomial complexity. The method requires a pseudo-concave lower-bound of \eqref{Prob:aSEE_ARIS_power_PCSI}, which would allow the use of fractional programming tools. To this end, we rewrite the SEE as follows
\begin{align}\label{Prob:SEE_ARIS_power_2}
	\text{SEE}(\boldsymbol{p}) &= g_{1, B}(\boldsymbol{p}) - g_{2, B}(\boldsymbol{p}) - g_{1, E}(\boldsymbol{p}) + g_{2, E}(\boldsymbol{p})\notag\\	
	&= g_{1,B}(\boldsymbol{p})+ g_{2,E}(\boldsymbol{p}) -  g_{2, B}(\boldsymbol{p}) - g_{1, E }(\boldsymbol{p}) 
\end{align}
wherein we define
\begin{align}
	&g_{1,B}(\boldsymbol{p})  = \frac{\sum_{k=1}^{K}\log_{2}\left(\!d_{k}^{(B)}+\sum_{k=1}^{K}p_{k}a^{(B)}_{k,k}\right)}{\sum_{k=1}^{K}\mu_{k}p_{k}+P_{c,eq}}\notag\\ 
	&g_{2,B}(\boldsymbol{p})  = \frac{\sum_{k=1}^{K}\log_{2}\left(d_{k}^{(B)}+\sum_{m\neq k}p_{m}a^{(B)}_{k,m}\right)}{\sum_{k=1}^{K}\mu_{k}p_{k}+P_{c,eq}}\notag
\end{align}
\begin{align}
    &g_{1,E}(\boldsymbol{p})  = \frac{\sum_{k=1}^{K}\log_{2}\left(\!d^{(E)}+\sum_{k=1}^{K}p_{k}a^{(E)}_{k}\right)}{\sum_{k=1}^{K}\mu_{k}p_{k}+P_{c,eq}}\notag\\ 
	&g_{2,E}(\boldsymbol{p})  = \frac{\sum_{k=1}^{K}\log_{2}\left(d^{(E)}+\sum_{m\neq k}p_{m}a^{(E)}_{m}\right)}{\sum_{k=1}^{K}\mu_{k}p_{k}+P_{c,eq}}\notag\;.
\end{align}
Then, define $f_{1,i}(\boldsymbol{p})$ and $f_{2,i}(\boldsymbol{p})$ as the numerators of $g_{1,i}(\boldsymbol{p})$ and $g_{2,i}(\boldsymbol{p})$, with $i=B,E$. Notably, while $f_{1,B}(\boldsymbol{p})$ and $f_{2,E}(\boldsymbol{p})$ are concave in $\boldsymbol{p}$, $-f_{2,B}(\boldsymbol{p})$ and $-f_{1,E}(\boldsymbol{p})$ render the numerator of the SEE non-concave.
Nevertheless, we can derive a pseudo-concave lower-bound for $\text{SEE}(\boldsymbol{p})$, denoted as $\widetilde{\text{SEE}}(\boldsymbol{p})$, by replacing $f_{2, B}(\boldsymbol{p})$ and $f_{1, E}(\boldsymbol{p})$ with their first-order Taylor expansion around any feasible point $\bar{\boldsymbol{p}}$. Thus,~\eqref{Prob:SEE_ARIS_power_2} becomes
\begin{align}
	&\text{SEE}(\boldsymbol{p}) \geq g_{1,B}(\boldsymbol{p})+ g_{2,E}(\boldsymbol{p}) -  g_{2, B}(\bar{\boldsymbol{p}}) - g_{1, E}(\bar{\boldsymbol{p}})\notag\\ 
	&-  \frac{\left(\nabla f_{2,B}(\boldsymbol{p})\right)^{T}\left(\boldsymbol{p}-\bar{\boldsymbol{p}}\right)}{\sum_{k=1}^{K}\mu_{k}p_{k}+P_{c,eq}} - \frac{\left(\nabla f_{1,E}(\boldsymbol{p})\right)^{T}\left(\boldsymbol{p}-\bar{\boldsymbol{p}}\right)}{\sum_{k=1}^{K}\mu_{k}p_{k}+P_{c,eq}} = \widetilde{\text{SEE}}(\boldsymbol{p})
\end{align}
wherein for all $j=1,\ldots,K$, it holds
\begin{align}
	\frac{\partial f_{2,B}}{\partial p_{j}}&=\sum_{k\neq j}\frac{a^{(B)}_{k,j}}{d_{k}^{(B)}+\sum_{m\neq k}p_{m}a^{(B)}_{k,m}}\\
	\frac{\partial f_{1,E}}{\partial p_{j}}&=\sum_{k = 1}^{K}\frac{a^{(B)}_{k,j}}{d^{(E)}+\sum_{m=1}^{K}p_{m}a^{(E)}_{m}}
\end{align}
Thus, a sequential fractional programming algorithm can be devised to address \eqref{Prob:SEE_ARIS_power_PCSI}, in which each iteration solves the problem
\begin{subequations}\label{Prob:SEE_ARIS_power_approx_PCSI}
	\begin{align}
		&\ds\max_{\boldsymbol{p}}\,\widetilde{\text{SEE}}(\boldsymbol{p})\\
		&\;\text{s.t.}\;0\leq p_{k}\leq P_{max,k}\;,\forall\; k=1,\ldots,K.
	\end{align}
\end{subequations}
Problem~\eqref{Prob:SEE_ARIS_power_approx_PCSI} is a pseudo-concave fractional program whose objective has a concave numerator and a convex denominator. Thus, it can be globally and efficiently solved by standard fractional programming techniques \cite{ZapNow15}.

\begin{algorithm}
	\caption{Power optimization - Perfect CSI}
	\begin{algorithmic}\label{Alg:SCA_p1_PCSI}
		\STATE $\epsilon > 0$, $\boldsymbol{\bar{p}}$ \texttt{any feasible vector};
		\REPEAT
		\STATE \texttt{Solve} \eqref{Prob:SEE_ARIS_power_approx_PCSI} \texttt{and let} $\boldsymbol{p^{\star}}$ \texttt{be the solution};
		\STATE $T=|\text{SEE}(\boldsymbol{p^{\star}})-\text{SEE}(\bar{\boldsymbol{p}})|$; \STATE $\bar{\boldsymbol{p}}=\boldsymbol{p^{\star}}$;
		\UNTIL{$T<\epsilon$}
	\end{algorithmic}
\end{algorithm}
\begin{proposition}\label{Prop:SCA_p1_PCSI}
	Algorithm \ref{Alg:SCA_p1_PCSI} monotonically improves the value of the objective of Problem \eqref{Prob:SEE_ARIS_power_PCSI} and converges to a point fulfilling the KKT optimality conditions of \eqref{Prob:SEE_ARIS_power_PCSI}.
\end{proposition}
\begin{IEEEproof}
	The proof follows similarly as for Proposition \ref{Prop:SCA_gamma1_RIS}, showing that Algorithm~\ref{Alg:SCA_p1_PCSI} fulfills all assumptions of the sequential fractional programming framework \cite{SeqCvxProg78}. 
\end{IEEEproof}

\subsection[Optimization of the LMMSE filters]{Optimization of $\boldsymbol{C}$}
The optimization of the receive filters in $\boldsymbol{C}$ exclusively influences the numerator of the SEE. Furthermore, it can be independently decoupled across users, simplifying the process to maximize each user's rate. The well-established solution to this problem is the linear minimum mean squared error (MMSE) receiver, which, for the considered case, is given by 
\begin{equation}
\boldsymbol{c}_{k}=\sqrt{p}_{k}\left(\sum_{m\neq k}p_{m}\boldsymbol{A}_{m}\boldsymbol{\gamma}\boldsymbol{\gamma}^{H}\boldsymbol{A}_{m}^{H}+\boldsymbol{W}\right)^{-1}\boldsymbol{A}_{k}\boldsymbol{\gamma}\;.
\end{equation}

\subsection[Overall Algorithm, Convergence and Complexity]{Overall  Algorithm, Convergence, Complexity, and practical implementation}
The overall alternating maximization algorithm can be stated as in Algorithm~\ref{Alg:SEE1}, which is guaranteed to converge, as shown in the coming Proposition \ref{Prop:SCA_Alt_PCSI}.
\begin{algorithm}[!h]
	\caption{Overall resource allocation - Perfect CSI}
	\begin{algorithmic}\label{Alg:SEE1}
		\STATE \texttt{Set} $\epsilon > 0$, \texttt{initialize} $\tilde{\boldsymbol{p}}$, $\tilde{\boldsymbol{\gamma}}$ \texttt{to 
			feasible values}
		\STATE \texttt{Compute} $\boldsymbol{c}_{k}=\sqrt{p}_{k}\boldsymbol{M}_{k}^{-1}\boldsymbol{A}_{k}$ \texttt{for all} $k$;
		\REPEAT
		\STATE \texttt{Compute} $\text{SEE}_{in}=\text{SEE}(\tilde{\boldsymbol{p}},\tilde{\boldsymbol{\gamma}},\boldsymbol{C})$;
		\STATE \texttt{Given} $\tilde{\boldsymbol{p}}$ \texttt{run} Algorithm \ref{Alg:SCA_gamma1_RIS};
		\STATE \texttt{Let} $\tilde{\boldsymbol{\gamma}}$ \texttt{be the optimized RIS vector};
		\STATE \texttt{Given} $\tilde{\boldsymbol{\gamma}}$ \texttt{run} Algorithm \ref{Alg:SCA_p1_PCSI};
		\STATE \texttt{Let} $\tilde{\boldsymbol{p}}$ \texttt{be the optimized power vector};
		\STATE \texttt{Compute} $\boldsymbol{c}_{k}=\sqrt{p}_{k}\boldsymbol{M}_{k}^{-1}\boldsymbol{A}_{k}$ \texttt{for all} $k$;
		\texttt{Compute} $\text{SEE}_{out}=\text{SEE}(\tilde{\boldsymbol{p}},\tilde{\boldsymbol{\gamma}},\boldsymbol{C})$
		\UNTIL{$|\text{SEE}_{out}-\text{SEE}_{in}|<\epsilon$}
	\end{algorithmic}
\end{algorithm}

\begin{proposition}\label{Prop:SCA_Alt_PCSI}
	Algorithm 3 monotonically increases the SEE value and converges in the value of the objective.
\end{proposition}

\begin{IEEEproof}
	Based on Propositions~\ref{Prop:SCA_gamma1_RIS} and~\ref{Prop:SCA_p1_PCSI}, we can infer that Algorithm~\ref{Alg:SEE1}  monotonically increases the SEE function in
	each step. Thus, since the SEE function has a finite maximizer, Algorithm~\ref{Alg:SEE1} eventually converges in the value of the objective.
\end{IEEEproof}

\begin{remark}\label{Rem:Special_PCSI}
Algorithm \ref{Alg:SEE1} can be specialized to maximize the SSR of the system. Indeed, the SSR coincides with the numerator of \eqref{Prob:aSEE_ARIS}, and, thus, it can be maximized by employing Algorithm \ref{Alg:SEE1} on a special instance of Problem \eqref{Prob:SEE_ARIS}, in which the denominator of \eqref{Prob:aSEE_ARIS} is replaced by $1$. 
\end{remark}

\textbf{Computational Complexity:}  Neglecting the complexity due to the computation of the closed-form receive filters in $\boldsymbol{C}$, the complexity of Algorithm~\ref{Alg:SEE1} can be obtained by recalling that a pseudo-concave maximization with $n$ variables can be restated as a concave maximization with $n+1$ variables \cite{ZapNow15}. Thus, the complexity of a pseudo-concave maximization problem with $n$ variables is polynomial in $n+1$. As a result, RIS and power optimization have complexity $\left(N+1\right)^\alpha$ and $\left(K+1\right)^\beta$, respectively \footnote{The exponents $\alpha$ and $\beta$ are not known, but for generic convex problems they can be bounded between 1 and 4~\cite{BenTal2001ConvexOptimization}}. Thus, the complexity of Algorithm~\ref{Alg:SEE1} is $\mathcal{C}_{1}=\mathcal{O}\left(I_{1}\left(I_{\gamma,1}\left(N+1\right)^\alpha + I_{p,1}\left(K+1\right)^\beta\right)\right)$ with $I_{\gamma,1}$, $I_{p,1}$ and $I_{1}$ the number of iteration that are required for Algorithms~\ref{Alg:SCA_gamma1_RIS},~\ref{Alg:SCA_p1_PCSI} and~\ref{Alg:SEE1} to reach convergence, respectively.

\textbf{Practical implementation:} Algorithm~\ref{Alg:SEE1} is meant to be implemented centrally by the BS, which then sends the configuration signal to the RIS and the transmit powers to be used to the mobile users'. The former requires the transmission of an amount of control information that scales with the number $N$ of RIS elements, while the latter requires the use of a control channel to transmit the $K$ power values. In addition, before Algorithm~\ref{Alg:SEE1} can be run, the BS needs to acquire the propagation channels. As for the channels that involve the legitimate nodes, they can be estimated by standard techniques for RIS-aided systems. Specifically, blind or pilot-based techniques can be used to estimate the product channels $\boldsymbol{A}_{k}$ for all $k=1,\ldots,K$ (see for example \cite{ZapTWC21,Buzzi21}), and the channel $\boldsymbol{G}_{B}$, which is a static channel since both RIS and BS are fixed nodes. From $\boldsymbol{A}_{k}$ and $\boldsymbol{G}_{B}$, the diagonal channel matrix $\boldsymbol{H}_{k}$ can be obtained since $\boldsymbol{A}_{k}=\boldsymbol{G}_{B}\boldsymbol{H}_{k}$. As for the channel $\boldsymbol{g}$ from the RIS to the eavesdropper, we recall that in the perfect CSI scenario, the eavesdropper is assumed to be a legitimate node, which, however, must not be allowed to decode the information signal. Thus, the channel $\boldsymbol{g}$ can also be estimated by standard methods. The overhead related to channel estimation depends heavily on the channel estimation routine that is used, but in general it will be proportional to the number of channel coefficients that must be estimated, thus scaling with $NN_{B}K$. 

\section{Optimization with statistical CSI} \label{Sec:OPT_SCSI}
In order to tackle Problem \eqref{Prob:SEE_ARIS_SCSI}, it is necessary to evaluate the statistical expectation in the numerator of \eqref{Prob:aSEE_ARIS_SCSI}. Unfortunately, a closed-form expression of the term $\mathbb{E}_{\boldsymbol{\delta}}\left[\log_{2}(1+\text{SINR}_{k,E})\right]$ is not available. Thus, to simplify the problem, we approximate the objective by taking the expectation inside the logarithm, which yields
\begin{align}
	&\mathbb{E}_{\boldsymbol{\delta}}\left[\log_{2}(1+\text{SINR}_{k,E})\right]\label{Eq:ApproxErgRate}\\
	&=\mathbb{E}_{\boldsymbol{\delta}}\left[\log_{2}\left(\sigma_{E}^{2}\!+\!\sigma_{RIS}^{2}\boldsymbol{g}_{E}^{H}\boldsymbol{\Gamma}\boldsymbol{\Gamma}^{H}\boldsymbol{g}_{E}\!+\!\sum_{m=1}^{K}p_{m}\left|\boldsymbol{g}_{E}^{H}\boldsymbol{H}_{m}\boldsymbol{\gamma}\right|^{2}\right)\right]\notag\\
	&-\mathbb{E}_{\boldsymbol{\delta}}\!\left[\!\log_{2}\left(\sigma_{E}^{2}\!+\!\sigma_{RIS}^{2}\boldsymbol{g}_{E}^{H}\boldsymbol{\Gamma}\boldsymbol{\Gamma}^{H}\boldsymbol{g}_{E}\!+\!\sum_{m\neq k}p_{m}\left|\boldsymbol{g}_{E}^{H}\boldsymbol{H}_{m}\boldsymbol{\gamma}\right|^{2}\!\right)\!\right]\notag\\
	&\approx \log_{2}\left(\sigma_{E}^{2}\!+\!\mathbb{E}_{\boldsymbol{\delta}}\!\left[\!\sigma_{RIS}^{2}\boldsymbol{g}_{E}^{H}\boldsymbol{\Gamma}\boldsymbol{\Gamma}^{H}\boldsymbol{g}_{E}\!+\!\sum_{m=1}^{K}p_{m}\left|\boldsymbol{g}_{E}^{H}\boldsymbol{H}_{m}\boldsymbol{\gamma}\right|^{2}\!\right]\!\right)\notag\\
	&-\log_{2}\left(\sigma_{E}^{2}\!+\!\mathbb{E}_{\boldsymbol{\delta}}\!\left[\!\sigma_{RIS}^{2}\boldsymbol{g}_{E}^{H}\boldsymbol{\Gamma}\boldsymbol{\Gamma}^{H}\boldsymbol{g}_{E}\!+\!\sum_{m\neq k}p_{m}\left|\boldsymbol{g}_{E}^{H}\boldsymbol{H}_{m}\boldsymbol{\gamma}\right|^{2}\!\right]\!\right)\notag
\end{align}
Since $\boldsymbol{g}_{E}\!=\!\boldsymbol{\widehat{g}}_{E}\!+\!\boldsymbol{\delta}$, elaborating, we obtain, for all $m=1,\ldots,K$
\begin{align}
	&\mathbb{E}_{\boldsymbol{\Delta}}\left[\left|\boldsymbol{g}_{E}^{H}\boldsymbol{H}_{m}\boldsymbol{\gamma}\right|^{2}\right]=\boldsymbol{\widehat{g}}_{E}^{H}\boldsymbol{H}_{m}\boldsymbol{\gamma}\boldsymbol{\gamma}^{H}\boldsymbol{H}_{m}^{H}\boldsymbol{\widehat{g}}_{E}+\sigma_{g}^{2}\|\boldsymbol{H}_{m}\boldsymbol{\gamma}\|^{2}\notag\\
	&=\boldsymbol{\gamma}^{H}\boldsymbol{H}_{m}^{H}\left(\boldsymbol{\widehat{g}}_{E}\boldsymbol{\widehat{g}}_{E}^{H}+\sigma_{g}^{2}\boldsymbol{I}_{N}\right)\boldsymbol{H}_{m}\boldsymbol{\gamma}=\|\boldsymbol{R}_{E}^{1/2}\boldsymbol{H}_{m}\gamma\|^{2}\\
	&\mathbb{E}_{\boldsymbol{\Delta}}\left[\boldsymbol{g}_{E}^{H}\boldsymbol{\Gamma}\boldsymbol{\Gamma}^{H}\boldsymbol{g}_{E}\right]=\boldsymbol{\gamma}^{H}\left(\boldsymbol{\widehat{g}}_{E}\boldsymbol{\widehat{g}}_{E}^{H}+\sigma_{g}^{2}\boldsymbol{I}_{N}\right)\boldsymbol{\gamma}\!=\!\|\boldsymbol{R}_{E}^{1/2}\boldsymbol{\gamma}\|^{2}, 
\end{align}
wherein $\boldsymbol{R}_{E}=\boldsymbol{\widehat{g}}_{E}\boldsymbol{\widehat{g}}_{E}^{H}+\sigma_{g}^{2}\boldsymbol{I}_{N}$.
Then, \eqref{Eq:ApproxErgRate} becomes 
\begin{equation}
	\mathbb{E}_{\boldsymbol{\delta}}\left[\log_{2}(1+\text{SINR}_{k,E})\right]\approx \log_{2}\left(1+\widetilde{\text{SINR}}_{k,E}\right)\;,
\end{equation}
with $\widetilde{\text{SINR}}_{k,E}$ given by 
\begin{equation}\label{Eq:SINRApprox}
	\hspace{-0.2cm}\widetilde{\text{SINR}}_{k,E}\!=\!\frac{p_{k}\|\boldsymbol{R}_{E}^{1/2}\boldsymbol{H}_{k}\boldsymbol{\gamma}\|^{2}}{\sum_{m\neq k}p_{m}\|\boldsymbol{R}_{E}^{1/2}\boldsymbol{H}_{m}\boldsymbol{\gamma}\|^{2}\!+\!\sigma_{RIS}^{2}\|\boldsymbol{R}_{E}^{1/2}\boldsymbol{\gamma}\|^{2}\!+\!\sigma_{E}^{2}}
\end{equation}
Thus, \eqref{Prob:aSEE_ARIS_SCSI} can be approximated as 
\begin{equation}
	\widetilde{\text{SEE}}=\frac{\sum_{k=1}^{K}\log_{2}(1+\text{SINR}_{k,B})-\log_{2}(1+\widetilde{\text{SINR}}_{k,E})}{\tr\left(\left(\boldsymbol{\gamma}\boldsymbol{\gamma}^{H}-\boldsymbol{I}_{N}\right)\boldsymbol{R}\right)+\sum_{k=1}^{K}p_{k}+P_{c}}
\end{equation}
and Problem \eqref{Prob:SEE_ARIS_SCSI} can be restated as
\begin{subequations}\label{Prob:SEE_ARIS_SCSI2}
	\begin{align}
		&\ds\max_{\boldsymbol{\gamma},\boldsymbol{p},\boldsymbol{C}}\; \widetilde{\text{SEE}}(\boldsymbol{\gamma},\boldsymbol{p},\boldsymbol{C})\label{Prob:aSEE_ARIS_SCSI2}\\
		&\;\text{s.t.}\;\tr\left(\boldsymbol{R}\right)\leq \tr(\boldsymbol{R}\boldsymbol{\gamma}\boldsymbol{\gamma}^{H})\leq P_{R,max}+\tr\left(\boldsymbol{R}\right)\label{Prob:bSEE_ARIS_SCSI2}\\
		&\;\quad\;\;0\leq p_{k}\leq P_{max,k}\;\forall\;k=1,\ldots,K\;,\label{Prob:cSEE_ARIS_SCSI2}
	\end{align}
\end{subequations}
Problem \eqref{Prob:SEE_ARIS_SCSI2} will be tackled again by alternating optimization of $\boldsymbol{\gamma}$, $\boldsymbol{p}$, and $\boldsymbol{C}$. 
\subsection{RIS optimization}
With respect to the RIS vector $\boldsymbol{\gamma}$, the problem is:
\begin{subequations}\label{Prob:SEC_ARIS_gamma}
	\begin{align}
	&\max_{\boldsymbol{\gamma}}\frac{\textstyle\sum_{k=1}^{K}\!\log_{2}\left(1\!+\!\text{SINR}_{k,B}\right) \!-\! \log_{2}\left(1\!+\!\widetilde{\text{SINR}}_{k,E}\right)}{\tr\left(\boldsymbol{R}\boldsymbol{\gamma}\boldsymbol{\gamma}^{H}\right)+P_{c,eq}} \label{Prob:aSEC_ARIS_gamma}\\
		&\;\text{s.t.}\;\tr\left(\boldsymbol{R}\right)\leq \tr(\boldsymbol{R}\boldsymbol{\gamma}\boldsymbol{\gamma}^{H})\leq P_{R,max}+\tr\left(\boldsymbol{R}\right)\label{Prob:bSEC_ARIS_gamma}
	\end{align}
\end{subequations}
wherein $P_{c,eq}=\sum_{k}\mu_{k}p_{k}+P_{c}-\tr(\boldsymbol{R})$. Problem \eqref{Prob:SEC_ARIS_gamma} presents similar difficulties as in the perfect CSI scenario, with \eqref{Prob:aSEC_ARIS_gamma} being 
not pseudo-concave and \eqref{Prob:bSEC_ARIS_gamma} being a non-convex constraint. These issues can be overcome by resorting to the sequential fractional programming method again, adopting similar lower bounds as in the perfect CSI scenario. To elaborate, let us recall that $\boldsymbol{c}_{k}^{H}\boldsymbol{W}\boldsymbol{c}_{k}=\sigma^{2}\|\boldsymbol{c}_{k}\|^{2}+\sigma_{RIS}^{2}\|\widetilde{\boldsymbol{U}}_{k}^{1/2}\boldsymbol{\gamma}\|^{2}$, and the definitions of $x_{B}$, $y_{B}$, $\bar{x}_{B}$ and $\bar{y}_{B}$ in \eqref{Eq:X_B}, \eqref{Eq:Y_B}, \eqref{Eq:X_B_bar}, and \eqref{Eq:Y_B_bar}. Next, we redefine the quantities $x_{E}$, $y_{E}$, $\bar{x}_{E}$, $\bar{y}_{E}$ as 
\begin{align}
	&x_{E}=p_{k}\|\boldsymbol{R}_{E}^{1/2}\boldsymbol{H}_{k}\boldsymbol{\gamma}\|^{2}\;,\; \bar{x}_{E}=p_{k}\|\boldsymbol{R}_{E}^{1/2}\boldsymbol{H}_{k}\boldsymbol{\bar{\gamma}}\|^{2}\\
	&y_{E}=\textstyle\sum_{m\neq k}p_{m}\|\boldsymbol{R}_{E}^{1/2}\boldsymbol{H}_{m}\boldsymbol{\gamma}\|^{2}+\sigma_{RIS}^{2}\|\boldsymbol{R}_{E}^{1/2}\boldsymbol{\gamma}\|^{2}\\
 &\bar{y}_{E}=\textstyle\sum_{m\neq k}p_{m}\|\boldsymbol{R}_{E}^{1/2}\boldsymbol{H}_{m}\boldsymbol{\bar{\gamma}}\|^{2}+\sigma_{RIS}^{2}\|\boldsymbol{R}_{E}^{1/2}\boldsymbol{\bar{\gamma}}\|^{2}\;.
\end{align}
Thus, for all $k=1,\ldots,K$, the $k$-th summand at the numerator of \eqref{Prob:aSEC_ARIS_gamma} can be written as in \eqref{Eq:SecRate_PCSI}, which, applying the bounds in \eqref{Eq:LB1} and \eqref{Eq:LB2}, can be lower-bounded, for any $\bar{\boldsymbol{\gamma}}$, as 
\begin{align}\label{Eq:LowerBound}
	&R_{s,k}^{SCSI}(\boldsymbol{\gamma})\geq \log_{2}\left(1+\frac{\bar{x}_{B}}{\bar{y}_{B}}\right)\!+\!\frac{\bar{x}_{B}}{\bar{y}_{B}}\left(\frac{2\sqrt{x_{B}}}{\sqrt{\bar{x}_{B}}}-\frac{x_{B}+y_{B}}{\bar{x}_{B}+\bar{y}_{B}}-1\right)\notag\\
	&\hspace{-0.3cm}+\log_{2}\left(1+\frac{\bar{y}_{E}}{\sigma_{E}^{2}}\right)\!+\!\frac{\bar{y}_{E}}{\sigma_{E}^{2}}\left(\frac{2\sqrt{y_{E}}}{\sqrt{\bar{y}_E}}-\frac{y_{E}+\sigma_{E}^{2}}{\bar{y}_{E}+\sigma_{E}^{2}}-1\right)\notag\\
	&\hspace{-0.3cm}+\log_{2}\left(\frac{\sigma_{E}^{2}}{\sigma_{E}^{2}+\bar{x}_{E}+\bar{y}_{E}}\right)\!-\!\frac{x_{E}+y_{E}-\bar{x}_{E}-\bar{y}_{E}}{\sigma_{E}^{2}+\bar{x}_{E}+\bar{y}_{E}}=\bar{R}_{s,k}^{SCSI}(\boldsymbol{\gamma})
\end{align}
All terms in \eqref{Eq:LowerBound} that depend on $\boldsymbol{\gamma}$ are concave, except for $\sqrt{x_{B}}$ and $\sqrt{y_{E}}$. However, as in the perfect CSI scenario, both $\sqrt{x_{B}}$ and $\sqrt{y_{E}}$ are convex and thus can be lower-bounded by their first-order Taylor expansion around $\boldsymbol{\bar{\gamma}}$. Specifically, $\sqrt{x_{B}}\geq \widetilde{x}_{B}$, with $\widetilde{x}_{B}$ given by \eqref{Eq:lowerXB}, while, as for $\sqrt{y_{E}}$, it can be written as 
\begin{align}
\sqrt{y_{E}}&=\sqrt{\boldsymbol{\gamma}^{H}\left(\sum_{m\neq k}p_{m}\boldsymbol{H}_{m}^{H}\boldsymbol{R}_{E}\boldsymbol{H_{m}}+\sigma_{RIS}^{2}\boldsymbol{R}_{E}\right)\boldsymbol{\gamma}}\notag\\
&=\left\|\left(\sum_{m\neq k}p_{m}\boldsymbol{H}_{m}^{H}\boldsymbol{R}_{E}\boldsymbol{H_{m}}+\sigma_{RIS}^{2}\boldsymbol{R}_{E}\right)^{1/2}\boldsymbol{\gamma}\right\|\;,
\end{align}
thus showing its convexity. Therefore, as in the perfect CSI scenario, it holds $\sqrt{y_{E}}\geq \widetilde{x}_{E}$ formally given by \eqref{Eq:lowerYE}, but with the matrix $\boldsymbol{M}_{k,E}$ redefined as 
\begin{equation}
\boldsymbol{M}_{k,E}=\left(\sum_{m\neq k}p_{m}\boldsymbol{H}_{m}^{H}\boldsymbol{R}_{E}\boldsymbol{H_{m}}+\sigma_{RIS}^{2}\boldsymbol{R}_{E}\right)\;.
\end{equation} 
Thus, each summand in the numerator of \eqref{Prob:aSEC_ARIS_gamma} can be lower-bounded by
\begin{align}
&R_{s,k}^{SCSI}(\boldsymbol{\gamma})\geq \widebar{R}_{s,k}^{SCSI}(\boldsymbol{\gamma}) \geq 
\log_{2}\left(1+\frac{\bar{x}_{B}}{\bar{y}_{B}}\right)\notag\\
&+\frac{\bar{x}_{B}}{\bar{y}_{B}}\left(\frac{2\widetilde{x}_{B}}{\sqrt{\bar{x}_{B}}}-\frac{x_{B}+y_{B}}{\bar{x}_{B}+\bar{y}_{B}}-1\right)+\log_{2}\left(1+\frac{\bar{y}_{E}}{\sigma_{E}^{2}}\right)\notag\\
&+\frac{\bar{y}_{E}}{\sigma_{E}^{2}}\left(\frac{2\widetilde{y}_{E}}{\sqrt{\bar{y}_E}}-\frac{y_{E}+\sigma_{E}^{2}}{\bar{y}_{E}+\sigma_{E}^{2}}-1\right)+\log_{2}\left(\frac{\sigma_{E}^{2}}{\sigma_{E}^{2}+\bar{x}_{E}+\bar{y}_{E}}\right)\notag\\
&-\frac{x_{E}+y_{E}-\bar{x}_{E}-\bar{y}_{E}}{\sigma_{E}^{2}+\bar{x}_{E}+\bar{y}_{E}}=\widetilde{R}_{s,k}^{SCSI}(\boldsymbol{\gamma})\;,
\end{align}
with the quantities $x_{E}$, $y_{E}$, $\bar{x}_{E}$, $\bar{y}_{E}$, $\boldsymbol{M}_{k,E}$ redefined as in this section. Finally, as for the left-hand side of \eqref{Prob:bSEC_ARIS_gamma}, it can be dealt with as in the perfect CSI scenario, which leads to the same expression as in \eqref{Eq:ConstraintSFP}. 
Then, in each iteration of the sequential method, the following problem must be globally solved:
\begin{subequations}\label{Prob:SEE_Gamma_Active}
	\begin{align}
		&\ds\max_{\boldsymbol{\gamma}}\;\frac{\textstyle\sum_{k=1}^{K}\widetilde{R}_{s,k}^{SCSI}}{\tr\left(\boldsymbol{R}\boldsymbol{\gamma}\boldsymbol{\gamma}^{H}\right)+P_{c,eq}} \label{Prob:aSEE_Gamma_Active}\\
		&\;\text{s.t.}\;\boldsymbol{\gamma}\boldsymbol{R}\boldsymbol{\gamma}^{H} \leq P_{R,max}+\tr\left(\boldsymbol{R}\right)\label{Prob:bSEE_Gamma_Active}\\
		&\; \quad \;\bar{\boldsymbol{\gamma}}\boldsymbol{R}\bar{\boldsymbol{\gamma}}^{H} + 2\Re\left\{\bar{\boldsymbol{\gamma}}^{H}\boldsymbol{R}\left(\boldsymbol{\gamma} - \bar{\boldsymbol{\gamma}}\right)\right\} \geq \tr\left(\boldsymbol{R}\right)\label{Prob:cSEE_Gamma_Active}\;,
	\end{align}
\end{subequations}
which can be efficiently accomplished by fractional programming techniques. The resulting procedure is given in Algorithm~\ref{Alg:SCA_gamma1}.

\begin{algorithm}
	\caption{RIS optimization - Statistical CSI}
	\begin{algorithmic}\label{Alg:SCA_gamma1}
		\STATE $\epsilon > 0$, $\boldsymbol{\bar{\gamma}}$ \texttt{any feasible vector};
		\REPEAT
        \STATE \texttt{Solve} \eqref{Prob:SEE_Gamma_Active} \texttt{and let} $\boldsymbol{\gamma^{\star}}$ \texttt{be the solution};
        \STATE $T=|\widetilde{\text{SEE}}(\boldsymbol{\gamma}^{\star})-\widetilde{\text{SEE}}(\bar{\boldsymbol{\gamma}})|$; 
        \STATE $\bar{\boldsymbol{\gamma}}=\boldsymbol{\gamma}^{\star}$;
		\UNTIL{$T<\epsilon$}
	\end{algorithmic}
\end{algorithm}

\begin{proposition}\label{Prop:SCA_gamma1}
	Algorithm \ref{Alg:SCA_gamma1} monotonically improves the value of the objective of Problem \eqref{Prob:SEC_ARIS_gamma} and converges to a KKT point of \eqref{Prob:SEC_ARIS_gamma}.
\end{proposition}
\begin{IEEEproof}
The proof follows along the same line of reasoning as in the perfect CSI scenario, since all assumptions of the sequential optimization framework are fulfilled. 
\end{IEEEproof}

\subsection{Transmit power optimization}
Let us define $b_{k}=\|\boldsymbol{R}_{E}^{1/2}\boldsymbol{H}_{k}\boldsymbol{\gamma}\|^{2}$, for all $k=1,\ldots,K$, and $q=\sigma_{RIS}^{2}\|\boldsymbol{R}_{E}^{1/2}\boldsymbol{\gamma}\|^{2}+\sigma_{E}^{2}$. Then, recalling the definitions of $a^{(B)}_{k,m}$, $d^{(B)}$, $\mu_{k}$, and $P_{c,eq}$ given in Section \ref{Sec:PowerPCSI}, the power optimization problem can be formulated as in \eqref{Prob:SEE_ARIS_power}, shown at the top of the next page.
\begin{figure*}
	\begin{equation}\label{Prob:SEE_ARIS_power}
		\max_{\{p_{k}\in[0,P_{max,k}]\}_{k=1}^{K}}\,\frac{\sum_{k=1}^{K}\ds\log_{2}\left(1+\frac{p_{k}a^{(B)}_{k,k}}{d_{k}^{(B)}+\sum_{m\neq k}p_{m}a^{(B)}_{k,m}}\right)-\log_{2}\left(1+\frac{p_{k}b_{k}}{q+\sum_{m\neq k}p_{m}b_{m}}\right)}{\sum_{k=1}^{K}\mu_{k} p_{k}+P_{c,eq}}
	\end{equation}
\end{figure*}
Problem \eqref{Prob:SEE_ARIS_power} can be addressed again by sequential fractional programming, observing that the SSR at the numerator of \eqref{Prob:SEE_ARIS_power} can be written as the difference of two concave functions of $\boldsymbol{p}$ as
\begin{align}
	\hspace{-0.2cm}R_{s}(\boldsymbol{p})&\!\!=\!\!\underbrace{\!\!\sum_{k=1}^{K}\log\!\left(\!d_{k}^{(B)}\!+\!\!\!\sum_{m=1}^{K}p_{m}a^{(B)}_{k,m}\right)\!+\!\log\!\left(\!q\!+\!\!\!\sum_{m\neq k}p_{m}b_{m}\!\!\right)}_{g_{1}(\boldsymbol{p})}\notag\\
	&\hspace{-1cm}\!-\!\underbrace{\!\left(\sum_{k=1}^{K}\!\!\log\!\left(\!d_{k}^{(B)}\!+\!\!\!\sum_{m\neq k}p_{m}a^{(B)}_{k,m}\right)\!\!+\!\log\!\left(\!q\!+\!\!\!\sum_{m=1}^{K}\!p_{m}b_{m}\right)\!\right)}_{g_{2}(\boldsymbol{p})}
\end{align}
Then a concave lower-bound of $R_{s}(\boldsymbol{p})$ can be obtained by linearizing $g_{2}(\boldsymbol{p})$ around any feasible point $\bar{\boldsymbol{p}}$, which yields the bound $R_{s}(\boldsymbol{p})\geq g_{1}(\boldsymbol{p})-g_{2}(\boldsymbol{\bar{p}})-\nabla g_{2}(\bar{\boldsymbol{p}}))^{T}(\boldsymbol{p}-\bar{\boldsymbol{p}})=\widetilde{R}_{s}(\boldsymbol{p})$, 
and thus, a surrogate problem that fits the assumptions of sequential fractional programming is obtained as
	\begin{align}\label{Prob:SEE_ARIS_power_approx}
		&\max_{\{p_{k}\in[0,P_{max,k}]\}_{k=1}^{K}}\,\frac{\widetilde{R}_{s}(\boldsymbol{p})}{\sum_{k=1}^{K}\mu_{k,eq}p_{k}+P_{c,eq}}
	\end{align}
Problem \eqref{Prob:SEE_ARIS_power_approx} can be efficiently solved by standard fractional programming methods, which leads to the power allocation subroutine stated in Algorithm \ref{Alg:SCA_p1}.

\begin{algorithm}
	\caption{Power optimization -  Statistical CSI}
	\begin{algorithmic}\label{Alg:SCA_p1}
		\STATE $\epsilon > 0$, $\boldsymbol{\bar{p}}$ \texttt{any feasible vector};
		\REPEAT
		\STATE \texttt{Solve} \eqref{Prob:SEE_ARIS_power_approx} \texttt{and let} $\boldsymbol{p^{\star}}$ \texttt{be the solution};
        \STATE $T=|\widetilde{\text{SEE}}(\boldsymbol{p}^{\star})-\widetilde{\text{SEE}}(\bar{\boldsymbol{p}})|$;
        \STATE $\bar{\boldsymbol{p}}=\boldsymbol{p}^{\star}$;
		\UNTIL{$T<\epsilon$}
	\end{algorithmic}
\end{algorithm}
\begin{proposition}\label{Prop:SCA_p1}
	Algorithm \ref{Alg:SCA_p1} monotonically improves the value of the objective of Problem \eqref{Prob:SEE_ARIS_power} and converges to a point fulfilling the KKT optimality conditions of \eqref{Prob:SEE_ARIS_power}.
\end{proposition}

\subsection{Receive filter optimization}
As in the perfect CSI scenario, the optimization of $\boldsymbol{C}$ exclusively influences the legitimate rate at the numerator of the SEE, and it can be decoupled across users, boiling down again to individual LMMSE filtering. 

\subsection[Overall Algorithm, Convergence and Complexity]{Overall  Algorithm, Convergence, Complexity, and Practical Implementation}
The overall maximization algorithm in the statistical CSI scenario is stated in Algorithm~\ref{Alg:SEE2}. 
\begin{algorithm}[!h]
	\caption{Overall resource allocation - Statistical CSI}
	\begin{algorithmic}\label{Alg:SEE2}
		\STATE \texttt{Set} $\epsilon > 0$, \texttt{initialize} $\tilde{\boldsymbol{p}}$, $\tilde{\boldsymbol{\gamma}}$ \texttt{to 
			feasible values}
		\STATE \texttt{Compute} $\boldsymbol{c}_{k}=\sqrt{p}_{k}\boldsymbol{M}_{k}^{-1}\boldsymbol{A}_{k}$ \texttt{for all} $k$;
		\REPEAT
		\STATE \texttt{Compute} $\widetilde{\text{SEE}}_{in}=\widetilde{\text{SEE}}(\tilde{\boldsymbol{p}},\tilde{\boldsymbol{\gamma}},\boldsymbol{C})$;
		\STATE \texttt{Given} $\tilde{\boldsymbol{p}}$ \texttt{run} Algorithm \ref{Alg:SCA_gamma1} \texttt{with output} $\tilde{\boldsymbol{\gamma}}$;
		\STATE \texttt{Given} $\tilde{\boldsymbol{\gamma}}$ \texttt{run} Algorithm \ref{Alg:SCA_p1} \texttt{with output} $\tilde{\boldsymbol{p}}$;
		\STATE \texttt{Compute} $\boldsymbol{c}_{k}\!\!=\!\!\sqrt{p}_{k}\boldsymbol{M}_{k}^{-1}\boldsymbol{A}_{k}$ and $\widetilde{\text{SEE}}_{out}\!\!=\!\!\widetilde{\text{SEE}}(\tilde{\boldsymbol{p}},\tilde{\boldsymbol{\gamma}},\boldsymbol{C})$;
		\UNTIL{$|\widetilde{\text{SEE}}_{out}-\widetilde{\text{SEE}}_{in}|<\epsilon$}
	\end{algorithmic}
\end{algorithm}

\begin{proposition}\label{Prop:SCA_Alt}
	Algorithm \ref{Alg:SEE2} monotonically increases the SEE value and converges in the value of the objective.
\end{proposition}

\begin{IEEEproof}
	Based on Propositions~\ref{Prop:SCA_gamma1} and~\ref{Prop:SCA_p1}, we can infer that Algorithm~\ref{Alg:SEE2} increases $\widetilde{\text{SEE}}$ in each step. Thus, since $\widetilde{\text{SEE}}$ has a finite maximizer, Algorithm~\ref{Alg:SEE2} eventually converges in the objective value.
\end{IEEEproof}

\begin{remark}\label{Rem:Special_SCSI}
Algorithm \ref{Alg:SEE2} can be specialized to maximize the ergodic (with respect to the channel $\boldmath{g}$) SSR of the system. Indeed, the ergodic SSR coincides with the numerator of \eqref{Prob:aSEE_ARIS_SCSI}, and thus it can be maximized by employing Algorithm \ref{Alg:SEE2} on a special instance of Problem \eqref{Prob:SEE_ARIS_SCSI}, in which the denominator of \eqref{Prob:aSEE_ARIS_SCSI} is replaced by $1$.  
\end{remark}

\textbf{Computational Complexity:}  Neglecting the complexity of computing the closed-form receive filters in $\boldsymbol{C}$, the complexity of Algorithm~\ref{Alg:SEE2} is obtained recalling that pseudo-concave maximizations with $n$ variables can be restated as concave maximization with $n+1$ variables and thus their complexity is polynomial in $n+1$~\cite{ZapNow15}. So, RIS and power optimization have complexity $\left(N+1\right)^\alpha$ and $\left(K+1\right)^\beta$,  respectively\footnote{The exponents $\alpha$ and $\beta$ are not known, but for generic convex problems they can be bounded between 1 and 4~\cite{BenTal2001ConvexOptimization}}, and, thus, the complexity of Algorithm~\ref{Alg:SEE2} is $\mathcal{C}_{1}=\mathcal{O}\left(I_{1}\left(I_{\gamma,1}\left(N+1\right)^\alpha + I_{p,1}\left(K+1\right)^\beta\right)\right)$ with  $I_{\gamma,1}$, $I_{p,1}$ and $I_{1}$ the number of iteration for Algorithms~\ref{Alg:SCA_gamma1},~\ref{Alg:SCA_p1},~\ref{Alg:SEE2} to converge.

\textbf{Practical implementation:} Similarly to Algorithm~\ref{Alg:SEE1}, Algorithm~\ref{Alg:SEE2} is also meant to be implemented centrally by the BS, and has similar control and channel estimation overhead requirements. The only major difference is that it allows for the possibility that a non-negligible error is made in the estimation of the channel $\boldsymbol{g}$ between the RIS and the eavesdropper. This scenario applies to situations in which the eavesdropper is a hidden node, whose position is approximately known. However, by letting $\sigma_{g}^{2}>>\|\boldsymbol{\widehat{g}}\|^{2}$, this case also applies to scenarios in which no information is available about the eavesdropper, which dispenses with the need of obtaining an estimate of $\boldsymbol{g}$.

\section{Numerical Analysis} \label{Sec:NUM_ANA}
Let us consider the setup detailed in Section~\ref{Sec:SysModel}, in which $K=4$, $N_{B}=4$, $N=100$, $B=20\,\textrm{MHz}$, $P_{0}=20\,\textrm{dBm}$, $P_{0,RIS}=30\,\textrm{dBm}$, $P_{c,n}=0\,\textrm{dBm}$. The noise spectral density is  $-174\,\textrm{dBm/Hz}$ with a noise figure of $5\,\textrm{dB}$, while the mobile users are distributed within a cell of radius of $50\,\textrm{m}$, with a minimum distance of $R_{n}=20\,\textrm{m}$ from the RIS and height from the ground randomly selected in $[1.5, 2.5]\,\textrm{m}$. The legitimate receiver is positioned $20\,\textrm{m}$ away from the RIS, while the eavesdropper is placed randomly within a $30\,\textrm{m}$ radius from the BS. Both the RIS and the BS are elevated at $10\,\textrm{m}$, while the eavesdropper has a random elevation in $[1.5, 2.5]\,\textrm{m}$. Thus, in this setup, the eavesdropper is considered another mobile node, which overhears the communication since it is in the coverage area. The power decay exponent of the channels from the mobile users to the RIS and from the eavesdropper to the RIS is $n_{h}=n_{g,E}=4$, while the power decay exponent of the channel from the RIS to the BS is $n_{g,B}=2$. All channels undergo Rician fading, with factor $K_{t}=4$ for the channel between the RIS and the legitimate receiver, and $K_{r}=2$ for the channels between the users and the eavesdropper to the RIS. This setup is motivated by the consideration that the RIS and BS are fixed and placed in favorable locations to ensure good propagation conditions. Indeed, the eavesdropper is either one of the network mobile nodes or, in any case, it is assumed to be in a less favorable propagation condition than the legitimate receiver. All results have been obtained by averaging over $10^3$ independent realizations of propagation channels and users' positions. For all illustrations regarding the performance with statistical CSI, we measure the quality of the channel estimation through the normalized error variance $\text{NEV} = \mathbb{E}[\|\boldsymbol{\delta}\|^2]/\mathbb{E}[\|\boldsymbol{g}\|^2]$. Unless otherwise specified, we set $\text{NEV}=0\,\textrm{dB}$.

\begin{figure}[!h]
	\centering
	\includegraphics[width=0.5\textwidth]{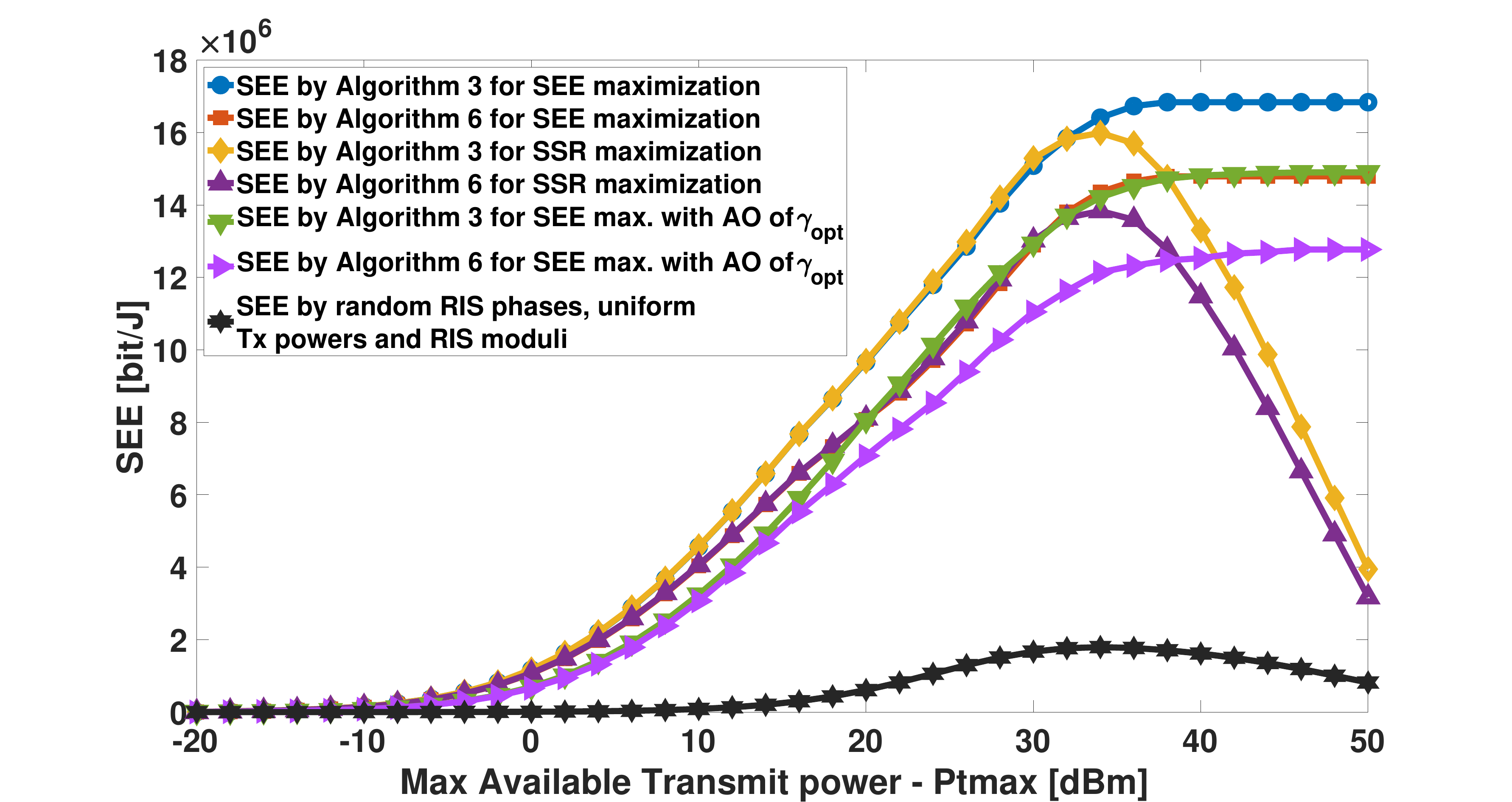}\caption{Achieved SEE versus $P_{tmax}$. $K=4$, $N_{B}=4$, $N=100$, $n_{h}=n_{g,E}=4, n_{g,B}=2$.} \label{fig:SEEvsP}
\end{figure}

Fig.~\ref{fig:SEEvsP} shows the SEE versus $P_{tmax}=P_{max,1}=\ldots=P_{max,K}$ for the following resource allocation schemes:
\begin{itemize}
\item[(a)] SEE in \eqref{Eq:SEE_PCSI} for the resource allocation obtained by running Algorithm \ref{Alg:SEE1} for SEE maximization.  
\item[(b)] SEE in \eqref{Eq:SEE_PCSI} for the resource allocation obtained by running Algorithm \ref{Alg:SEE2} for SEE maximization. Specifically, in this case, the SEE function in \eqref{Eq:SEE_PCSI} is computed for the true realizations of the channel $\boldsymbol{g}$, but using the RIS coefficients, transmit powers, and receive filters output by running Algorithm \ref{Alg:SEE2} for SEE maximization. Thus, statistical CSI is assumed at the design stage, but performance analysis is done for the true channel realization of $\boldsymbol{g}$. 
\item[(c)] SEE in \eqref{Eq:SEE_PCSI} for the resource allocation obtained by running Algorithm \ref{Alg:SEE1} for SSR maximization \footnote{Recall that Algorithm \ref{Alg:SEE1} can be specialized to maximize the SSR, as explained in Remark \ref{Rem:Special_PCSI}}. 
\item [(d)] SEE in \eqref{Eq:SEE_PCSI} for the resource allocation obtained by running Algorithm \ref{Alg:SEE2} for SSR maximization \footnote{Recall that Algorithm \ref{Alg:SEE2} can be specialized to maximize the ergodic SSR, as explained in Remark \ref{Rem:Special_SCSI}}. As for Scheme (b), also in this case, the SEE function in \eqref{Eq:SEE_PCSI} is computed for the true realizations of $\boldsymbol{g}$, but using the RIS coefficients, transmit powers, and receive filters output by running Algorithm \ref{Alg:SEE2} for ergodic SSR maximization.

\item[(e)] SEE in \eqref{Eq:SEE_PCSI} obtained by a modified version of Algorithm \ref{Alg:SEE1} for SSE maximization, in which the RIS reflection vector $\boldsymbol{\gamma}$ is optimized following the method from \cite{WuTCOM20}, i.e. by alternating optimization of its components. For each $n$, the modulus and phase of the component $\gamma_{n}$ have been optimized by exhaustive search.
\item[(f)] SEE in \eqref{Eq:SEE_PCSI} obtained by a modified version of Algorithm \ref{Alg:SEE2} for SSE maximization, in which the RIS reflection vector $\boldsymbol{\gamma}$ is optimized following the method from \cite{WuTCOM20}, i.e. by alternating optimization of its components. For each $n$, the modulus and phase of the component $\gamma_{n}$ have been optimized by exhaustive search.
\item[(g)] SEE in \eqref{Eq:SEE_PCSI} obtained by allocating the phases of the RIS coefficients randomly in $[0,2\pi]$, and employing a uniform allocation of the moduli of the RIS coefficients and the transmit powers. Specifically, $|\gamma_{n}|^{2}=1+P_{R,max}/\tr{(\boldsymbol{R})}$ for all $n$, and $p_{k}=P_{tmax}/K$ for all $k$. 
\end{itemize}
The results show that, as expected, the lack of perfect CSI of the eavesdropper's channel causes performance degradation for all resource allocation schemes. However, the gap is not severe, and confidential communications can be ensured with SEE value higher than $10\,\textrm{Mb/J}$ even when statistical CSI is assumed. Moreover, it is seen that the use of Schemes (a) and (b) allows saving a significant amount of energy, especially for high values of $P_{tmax}$, compared to Schemes (c) and (d) that pursue the maximization of the SSR. Finally, all schemes provide much higher SEE values than the heuristic Scheme (e) that does not perform any optimization.

\begin{figure}[!h]
	\centering
	\includegraphics[width=0.5\textwidth]{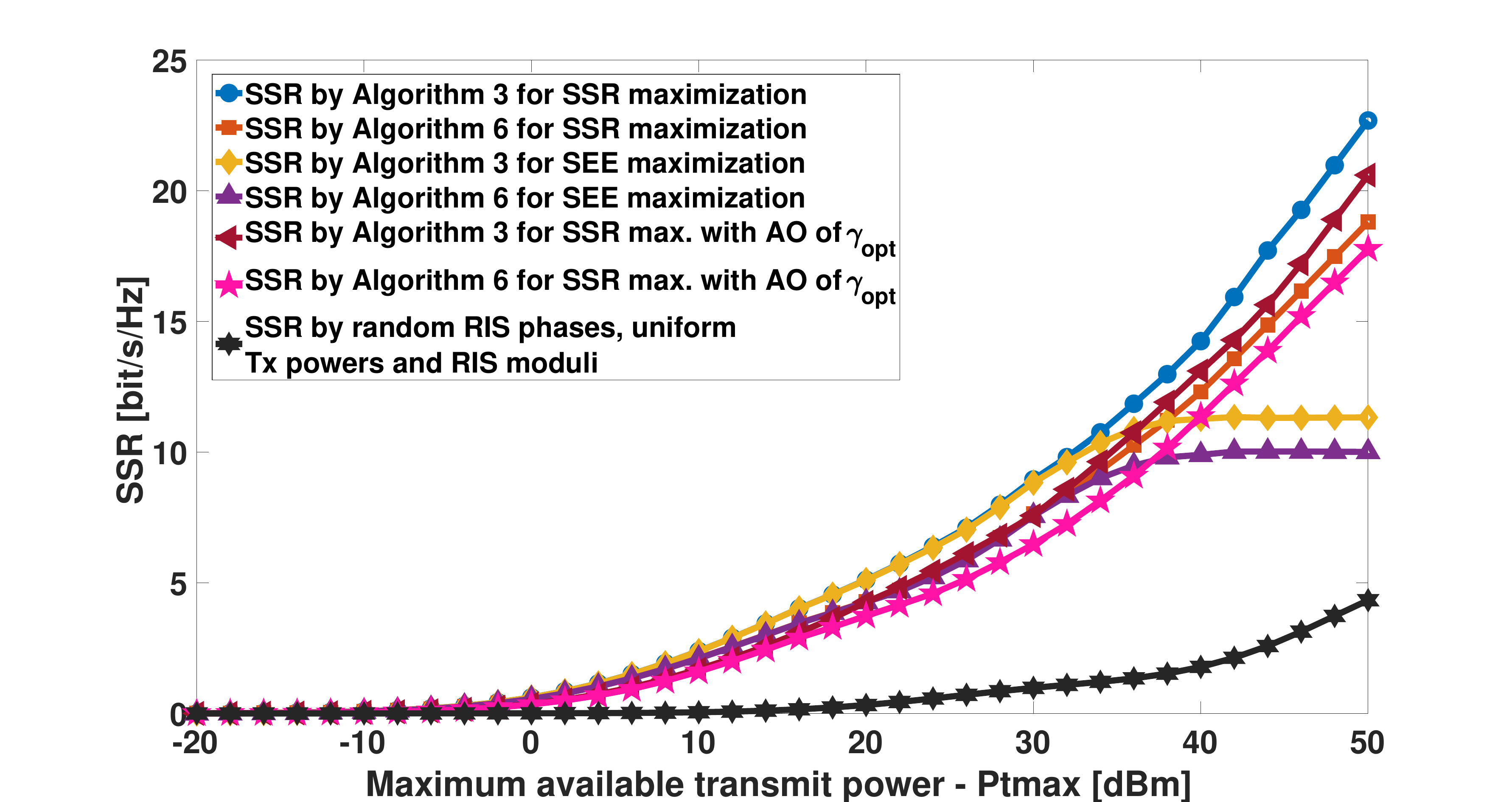}\caption{Achieved SSR versus $P_{tmax}$. $K=4$, $N_{B}=4$, $N=100$, $n_{h}=n_{g,E}=4, n_{g,B}=2$.} \label{fig:SSRvsP}
\end{figure}

Fig.~\ref{fig:SSRvsP} considers the same resource allocation schemes as in Fig. \ref{fig:SEEvsP}, with the difference that the SSR in \eqref{Eq:SecRate_PCSI} is shown instead of the SEE. Similar considerations as for Fig.~\ref{fig:SEEvsP} can be made. In particular, it can be seen that in the case of statistical CSI, the SSR saturates for the considered range of $P_{tmax}$, showing a limited gap between the maximization of the SSR and the maximization of the SEE. Thus, when statistical CSI is available, a very limited gap in terms of SSR is observed by allocating the system resources for SEE maximization. 

\begin{figure}[!h]
	\centering
	\includegraphics[width=0.5\textwidth]{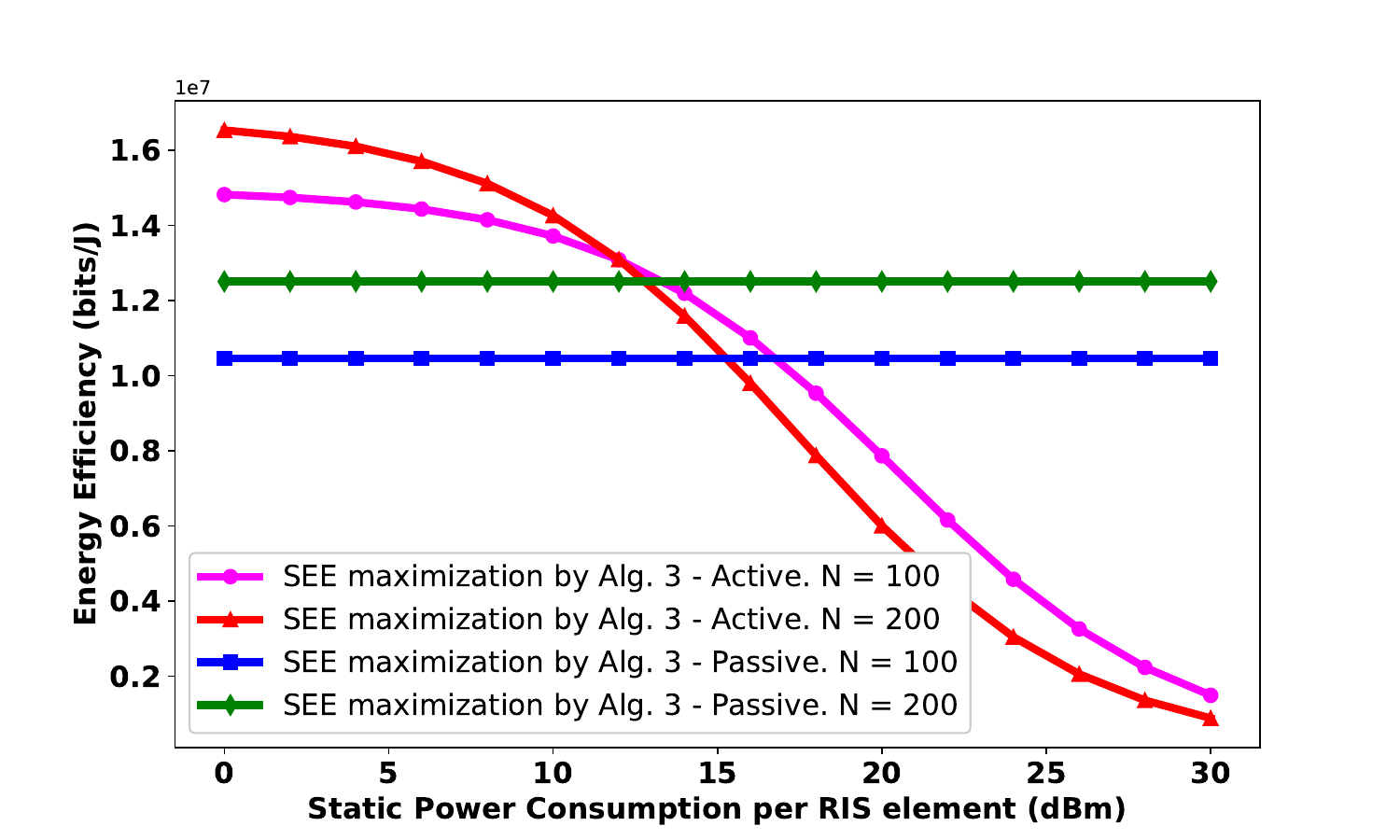}\caption{Achieved SEE for active and nearly-passive RIS versus static power consumption. $K=4$, $N_{B}=4$, $N=(100,200)$, $P_{tmax}=30\,\textrm{dBm}$, $P^{(p)}_{c,n}=0\,\textrm{dBm}$, $P^{(a)}_{0,RIS}=20\,\textrm{dBm}$, and $P^{(p)}_{0,RIS}=10\,\textrm{dBm}$.} \label{fig:EEvsPcn_pCSI}
\end{figure}

\begin{figure}[!h]
	\centering
	\includegraphics[width=0.5\textwidth]{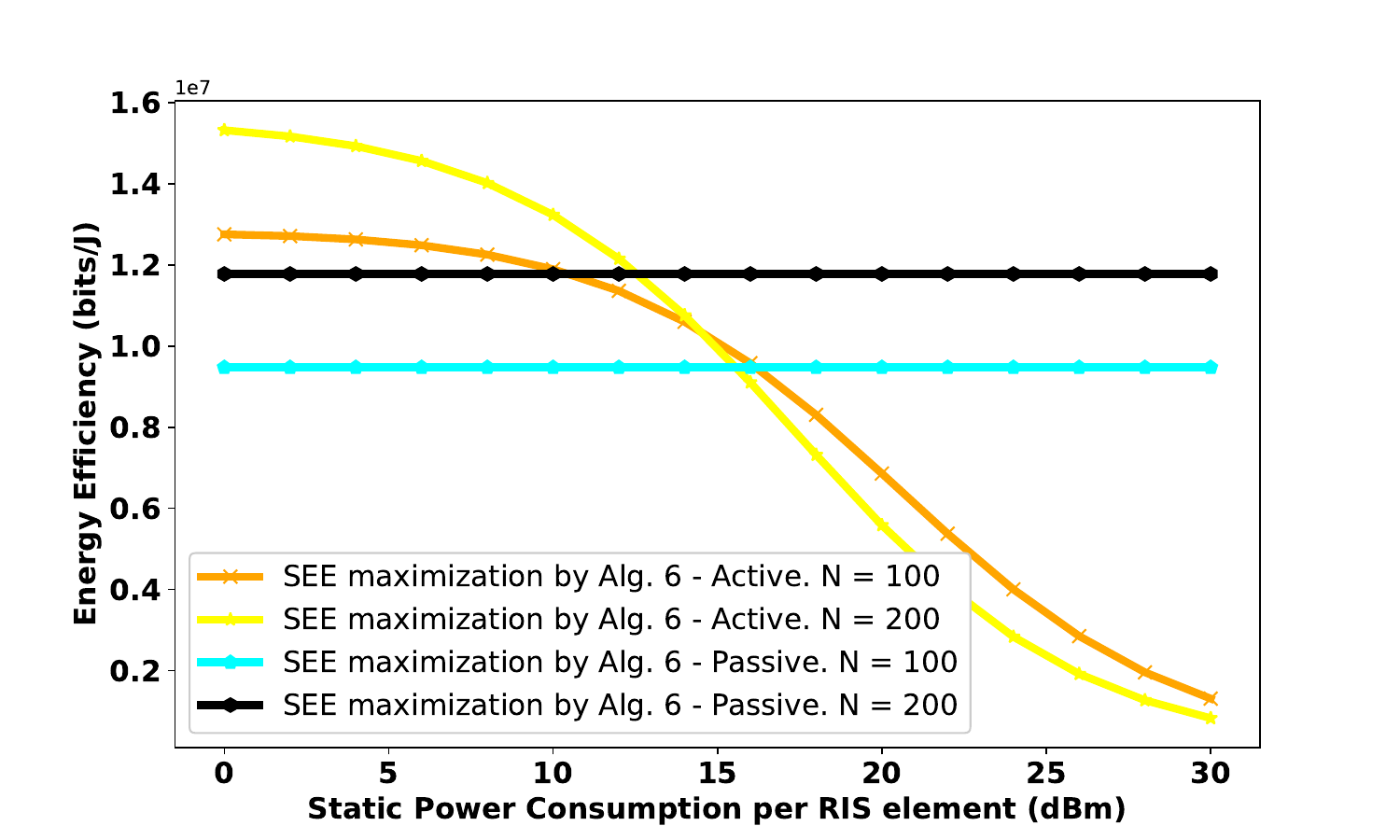}\caption{Achieved SEE for active and nealry-passive RIS versus static power consumption. $K=4$, $N_{B}=4$, $N=(100,200)$, $P_{tmax}=30\,\textrm{dBm}$, $P^{(p)}_{c,n}=0\,\textrm{dBm}$, $P^{(a)}_{0,RIS}=20\,\textrm{dBm}$, and $P^{(p)}_{0,RIS}=10\,\textrm{dBm}$.} \label{fig:EEvsPcn_sCSI}
\end{figure}

Figs.~\ref{fig:EEvsPcn_pCSI} and \ref{fig:EEvsPcn_sCSI} compare the SEE in \eqref{Eq:SEE_PCSI} achieved by using active and nearly-passive RISs. Specifically, Fig.~\ref{fig:EEvsPcn_pCSI} considers the SEE obtained with perfect CSI by Algorithm~\ref{Alg:SEE1}, while Fig.~\ref{fig:EEvsPcn_sCSI} considers the SEE obtained with statistical CSI by Algorithm~\ref{Alg:SEE2}. In both cases, it is assumed that the per-element static power consumption of the nearly-passive RIS is $P^{(p)}_{c,n}=0\,\textrm{dBm}$, while the per-element power consumption of the active RIS is assumed to vary from $P^{(a)}_{c,n}=0\,\textrm{dBm}$ to  $P^{(a)}_{c,n}=40\,\textrm{dBm}$, due to the presence of the analog amplifiers. Moreover, the rest of the RIS static power consumption is assumed to be $P_{0,RIS}^{(a)}=20\,\textrm{dBm}$ for the active RIS, and $P_{0,RIS}^{(a)}=10\,\textrm{dBm}$ for the nearly-passive RIS, always owing to the presence of the additional hardware for signal amplification in the active RIS. The comparison is made for $N=100$ and $N=200$. In both cases, as $P^{(a)}_{c,n}$ increases, the active RIS becomes less energy-efficient than the nearly-passive one, with the crossing point occurring for lower values of $P^{(a)}_{c,n}$ when $N$ is larger. Thus, as expected, there exists a trade-off between active and nearly-passive RIS as far as SEE is concerned, and thus, there exist practical system configurations for which nearly-passive RISs are to be preferred to achieve a higher SEE. 

\begin{figure}[!h]
	\centering
	\includegraphics[width=0.5\textwidth]{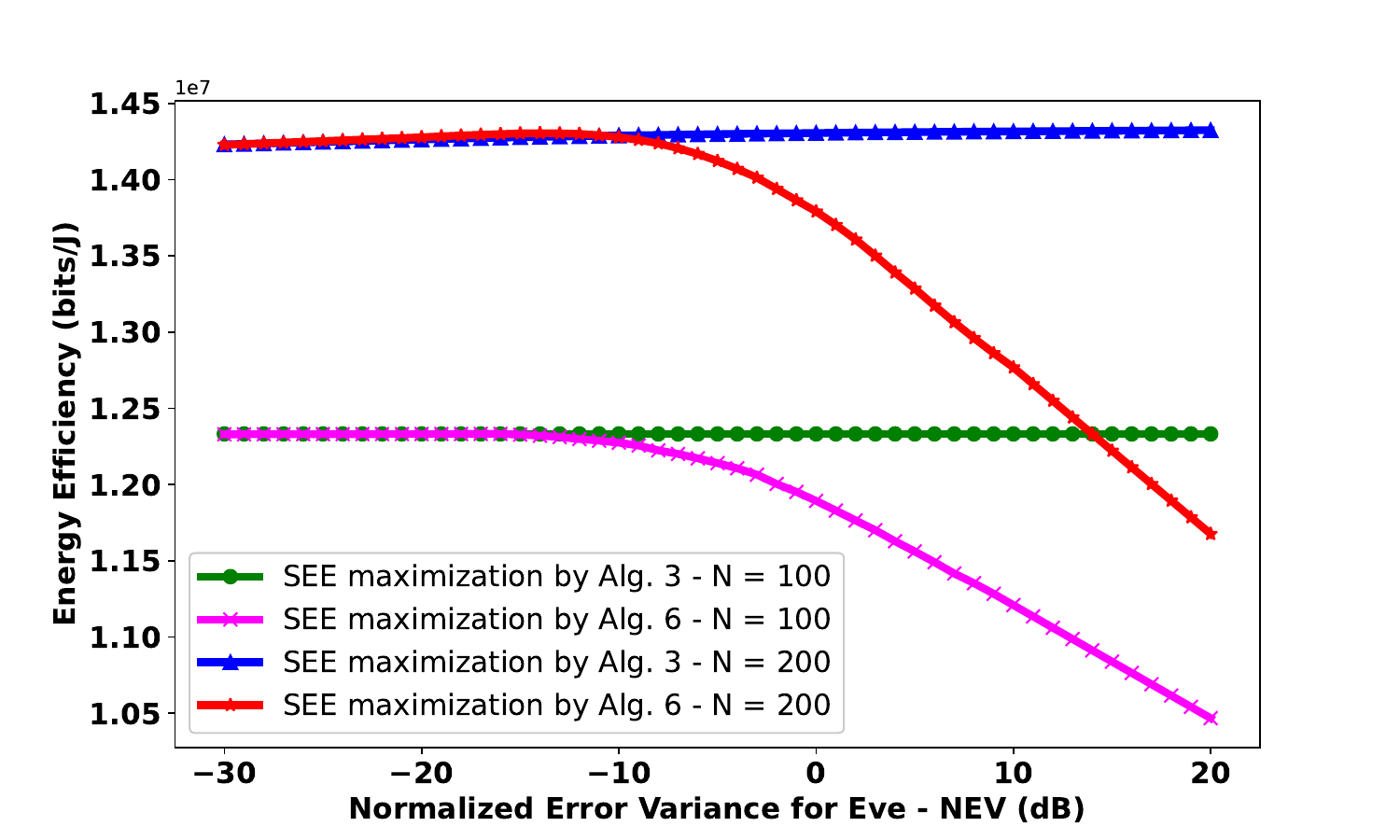}\caption{Achieved SEE versus NEV. $K=4$, $N_{B}=4$, $P_{tmax}=30\,\textrm{dBm}$, $N=100$, $n_{h}=n_{g,E}=4, n_{g,B}=2$.} \label{fig:SSRvsNEV}
\end{figure}

Fig. \ref{fig:SSRvsNEV} shows the SEE by Algorithm~\ref{Alg:SEE1} and Algorithm~\ref{Alg:SEE2} versus the NEV. The cases with $N=100$ and $N=200$ are considered for both scenarios. While the SEE obtained with perfect CSI does not depend on the NEV, the SEE obtained with statistical CSI decreases as the NEV increases. However, the gap remains negligible for NEV smaller than $0\,\textrm{dB}$. On the other hand, for larger values of the NEV, the performance obtained through Algorithm \ref{Alg:SEE2} deteriorates. 

\subsection{Relaxation and Projection}
In practice, RISs cannot provide any continuous value of reflection coefficients. This can be addressed by quantizing the phases and moduli of the coefficients output by Algorithm \ref{Alg:SEE1} and Algorithm~\ref{Alg:SEE2}. This approach is known as \emph{relaxation and projection}, because it solves the problem in the relaxed continuous domain, and then projects the solution in the discrete domain. The main advantage of this method is its low complexity compared to other approaches that directly optimize in the discrete domain. To elaborate on this, let us observe that the quantization step can be performed separately for each reflection coefficient and thus has a linear complexity in $N$. Thus, recalling that the complexity of Algorithms \ref{Alg:SEE1} and \ref{Alg:SEE2} has been shown to be polynomial in $N$, the complexity of the relaxation and projection method is polynomial in $N$, too. Instead, optimizing directly in the discrete domain would have an exponential complexity in the number of discrete variables $N$.

In practice, for all $n=1,\ldots,N$, the phase of $\gamma_{n}$ is quantized in $[0,2\pi]$. Instead, as for the modulus of $\gamma_{n}$, identifying the quantization range is less straightforward due to the constraint in \eqref{Prob:bSEE_ARIS}. Clearly, a simple approach would be to bound $|\gamma_{n}|^{2}\in[0,\tr({\boldsymbol{R}})+P_{R,max}]$. However, this yields a too-large range for $|\gamma_{n}|$, which leads to unsatisfactory performance. For this reason, we have adopted the following approach. First, we observe that $
R_{min}\|\boldsymbol{\gamma}\|^{2}\leq\tr({\boldsymbol{R}}\boldsymbol{\gamma}\boldsymbol{\gamma}^{H})\leq R_{max}\|\boldsymbol{\gamma}\|^{2}$, with $R_{max}$ and $R_{min}$ the maximum and minimum values of the diagonal elements of $\boldsymbol{R}$. Then, exploiting \eqref{Prob:bSEE_ARIS} yields
\begin{equation}
\frac{\tr({\boldsymbol{R}})}{R_{max}}\leq\|\boldsymbol{\gamma}\|^{2}\leq \frac{\tr({\boldsymbol{R}})+P_{Rmax}}{R_{min}}\;,
\end{equation}
from which, making the approximation  $\|\boldsymbol{\gamma}\|^{2}\approx N |\gamma_{n}|^{2}$, for all $n=1,\ldots,N$, we finally obtain that, for all $n$, 
\begin{equation}
|\gamma_{n}|\in\left[\sqrt{\frac{\tr({\boldsymbol{R}})}{R_{max}N}},\sqrt{\frac{\tr({\boldsymbol{R}})+P_{R,max}}{R_{min}N}}\right]\;,
\end{equation}
which defines the interval that has been used to quantize the output of Algorithms \ref{Alg:SEE1} and \ref{Alg:SEE2}. 

\begin{figure}[!h]
	\centering
	\includegraphics[width=0.5\textwidth]{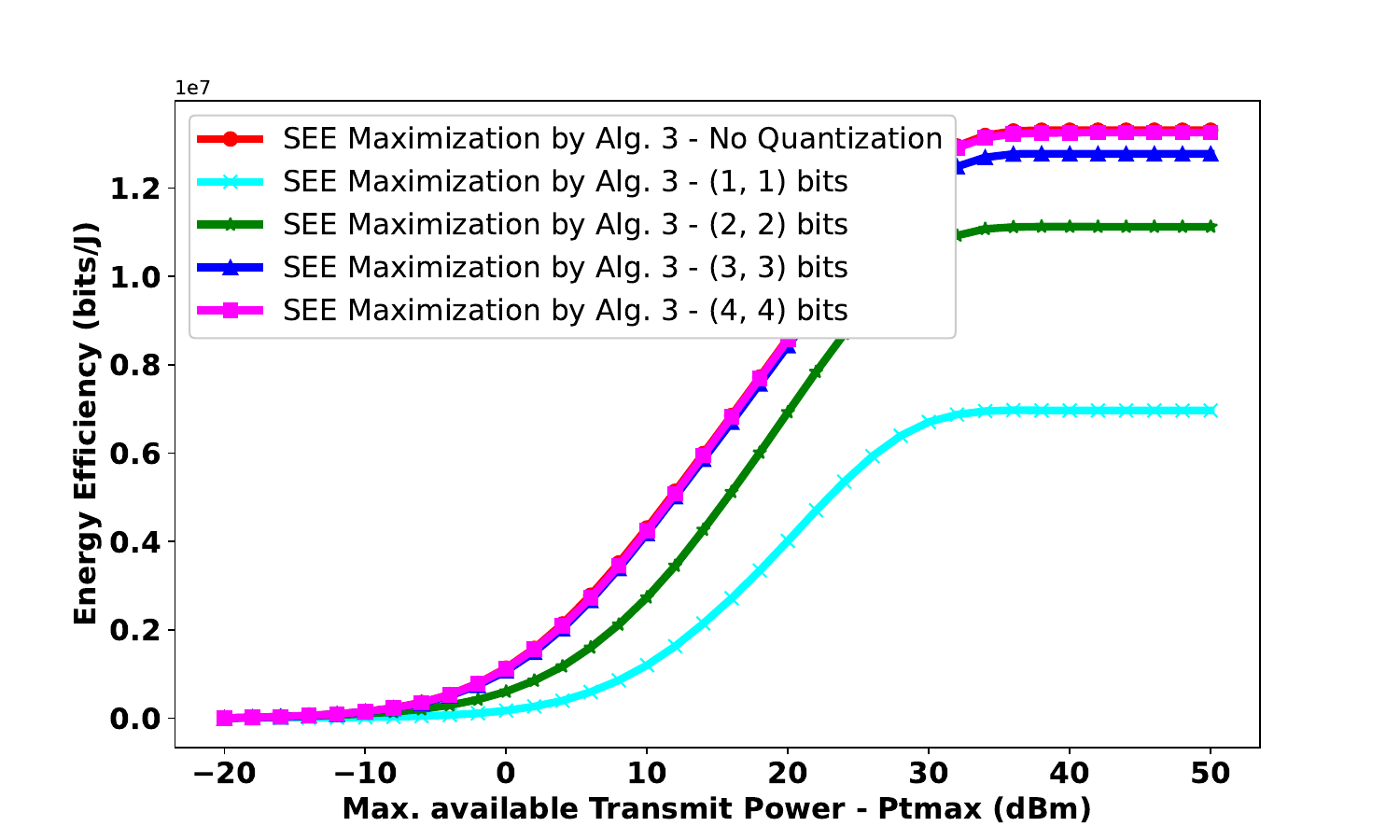}\caption{Achieved quantized SEE by Alg. \ref{Alg:SEE1} versus $P_{tmax}$ for different values of quantization bits (bit\_phase, bit\_amplitude). $K=4$, $N_{B}=4$, $N=100$, $n_{h}=n_{g,E}=4, n_{g,B}=2$.} \label{fig:SEEvsP_Q_pCSI}
\end{figure}

\begin{figure}[!h]
	\centering
	\includegraphics[width=0.5\textwidth]{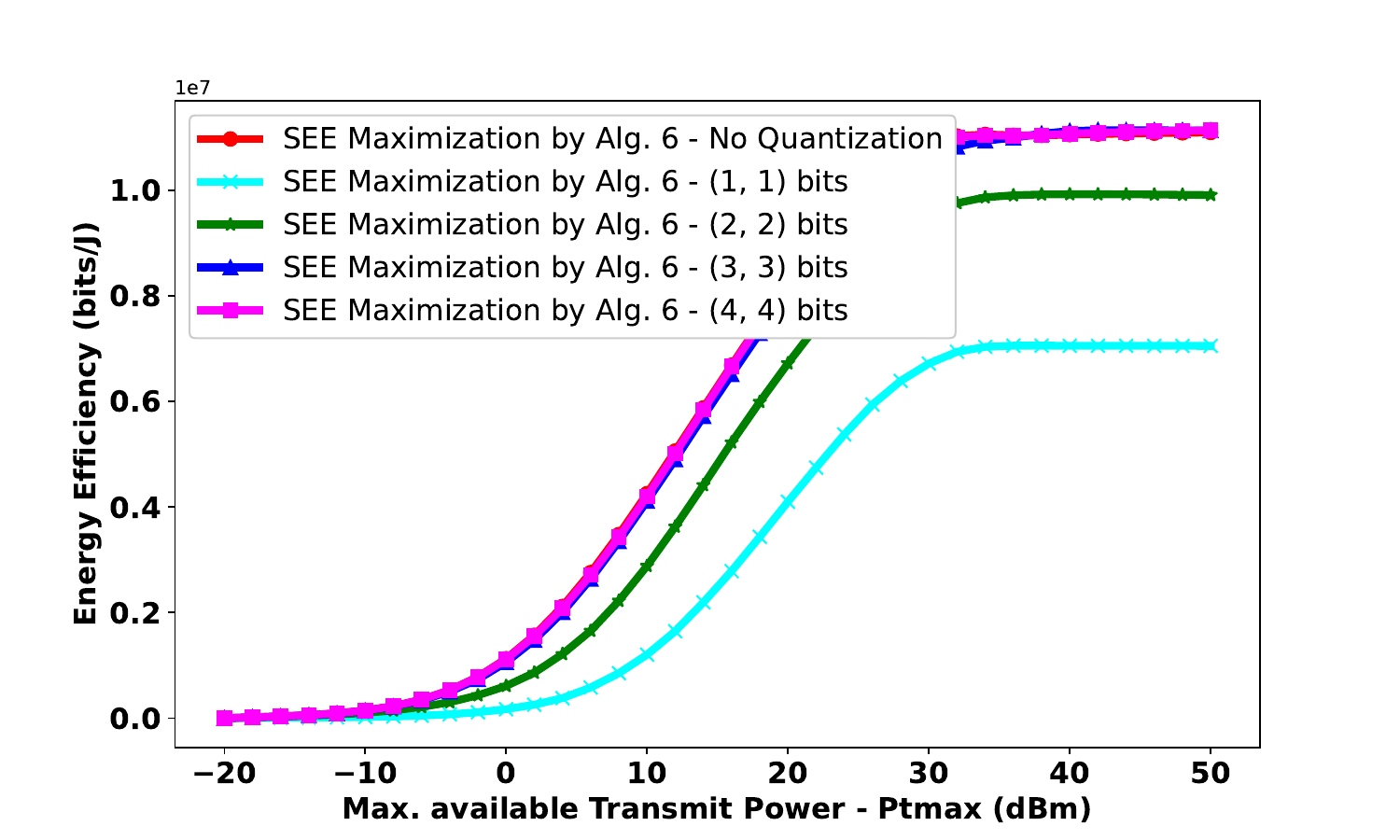}\caption{Achieved Quantized SEE by Alg. \ref{Alg:SEE2} versus $P_{tmax}$ for different values of quantization bits (bit\_phase, bit\_amplitude). $K=4$, $N_{B}=4$, $N=100$, $n_{h}=n_{g,E}=4, n_{g,B}=2$.} \label{fig:SEEvsP_Q_sCSI}
\end{figure}

Given this context, Figs. \ref{fig:SEEvsP_Q_pCSI} and \ref{fig:SEEvsP_Q_sCSI} address the performance of the relaxation and projection method, as far as the SEE is concerned. Specifically, Fig. \ref{fig:SEEvsP_Q_pCSI} considers perfect CSI and shows the SEE obtained by quantizing with $1$, $2$, $3$, and $4$ bits the output of Algorithm \ref{Alg:SEE1} for SEE maximization, while Fig. \ref{fig:SEEvsP_Q_sCSI} considers statistical CSI and shows the SEE obtained by quantizing with $1$, $2$, $3$, and $4$ bits the output of Algorithm \ref{Alg:SEE2} for SEE maximization. The performance obtained without quantization (i.e., employing the continuous output of the algorithms) is also reported, which provides an upper-bound to the performance of the relaxation and restriction method. The results show that, with both perfect and statistical CSI, using $3$ or more quantization bits leads to a negligible performance loss. Instead, the use of $2$ bits leads to a visible but still acceptable  loss, while using $1$ bit causes a significant SEE degradation. Moreover, it is seen that the gap is smaller in the statistical CSI regime, which can be explained observing that perfect allocation of the RIS coefficients is less relevant when perfect CSI of the RIS channel to the eavesdropper is not available. 

\begin{figure}[!h]
	\centering
	\includegraphics[width=0.5\textwidth]{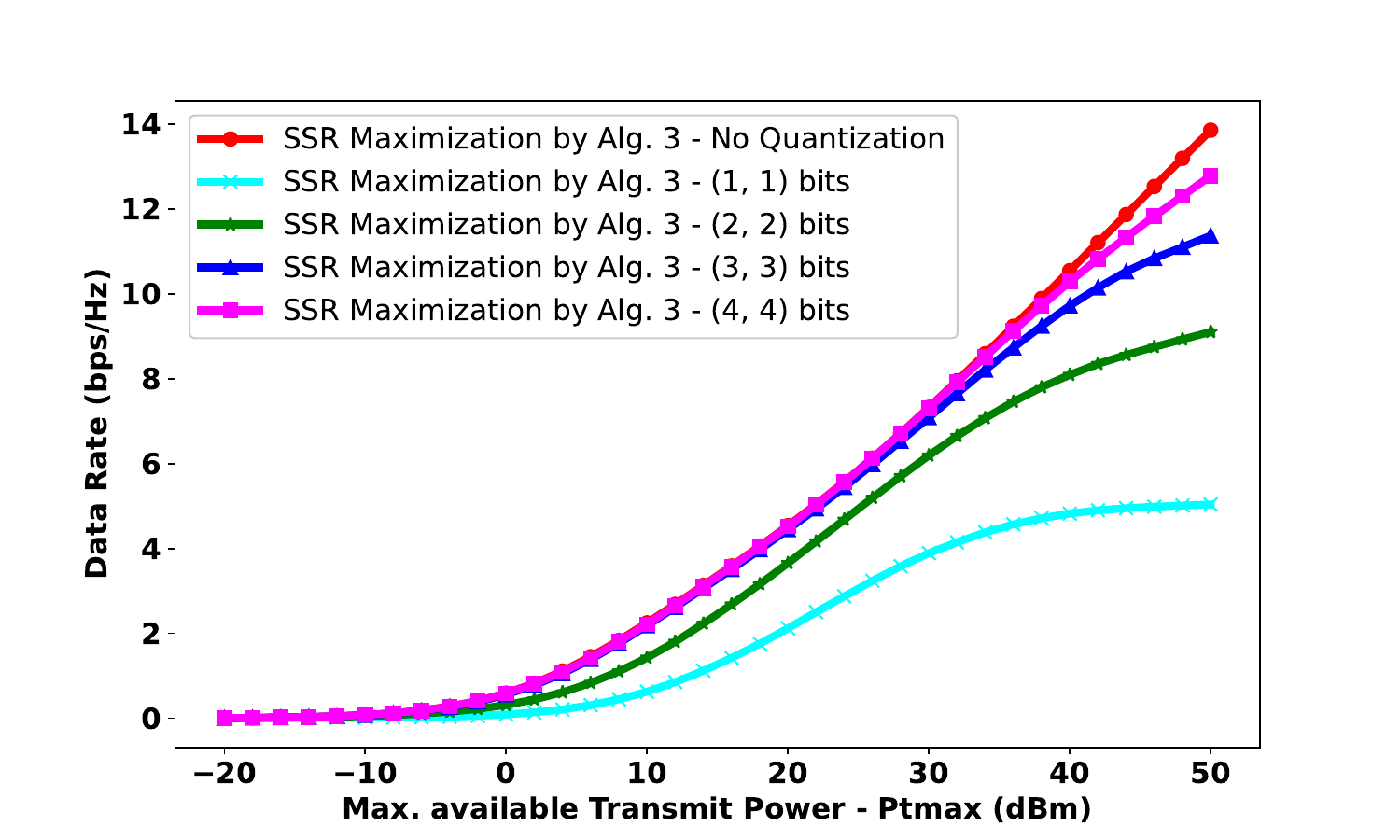}\caption{Achieved Quantized SSR by Alg. \ref{Alg:SEE1} versus $P_{tmax}$ for different values of quantization bits (bit\_phase, bit\_amplitude). $K=4$, $N_{B}=4$, $N=100$, $n_{h}=n_{g,E}=4, n_{g,B}=2$.} \label{fig:SSRvsP_Q_pCSI}
\end{figure}

\begin{figure}[!h]
	\centering
	\includegraphics[width=0.5\textwidth]{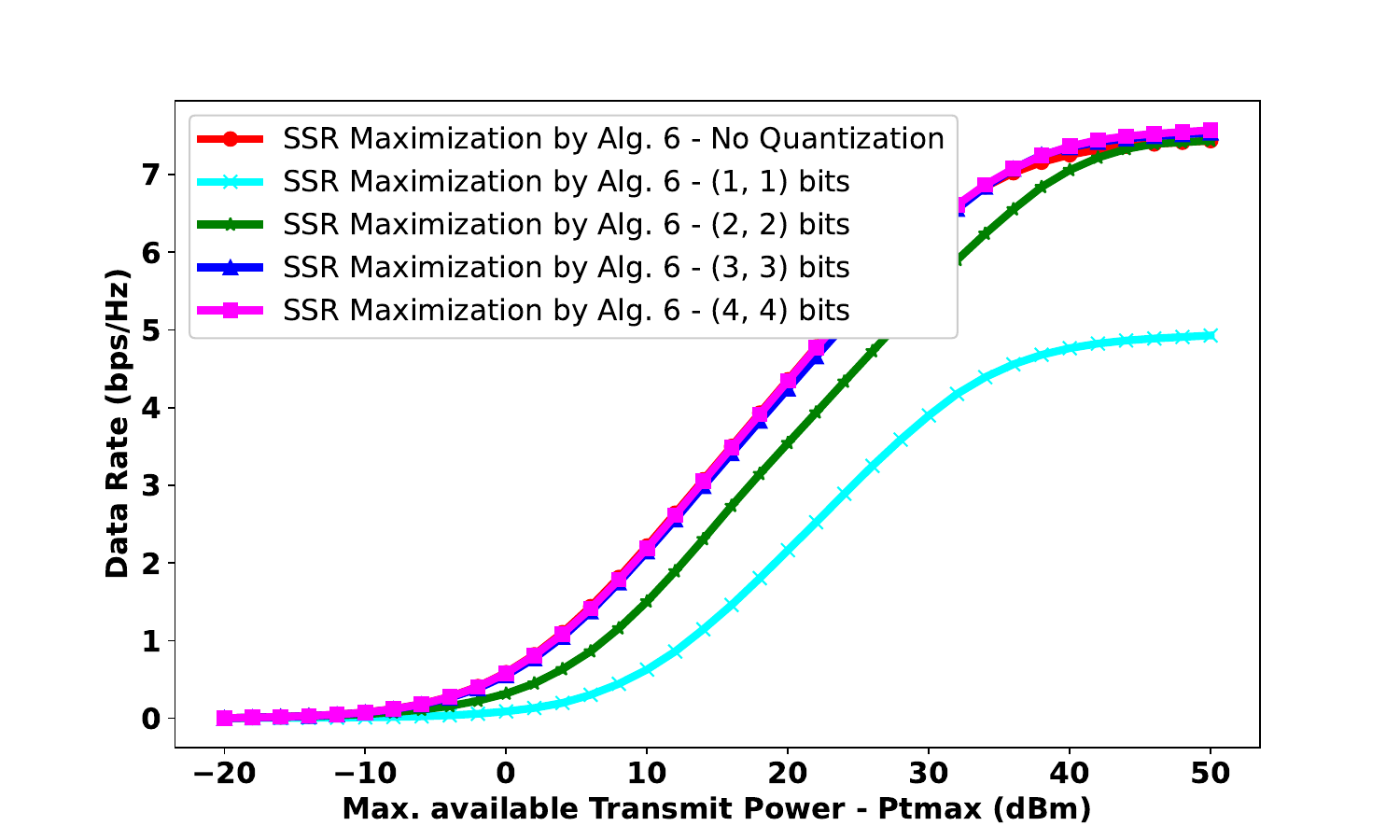}\caption{Achieved Quantized SSR by Alg. \ref{Alg:SEE2} versus $P_{tmax}$ for different values of quantization bits (bit\_phase, bit\_amplitude). $K=4$, $N_{B}=4$, $N=100$, $n_{h}=n_{g,E}=4, n_{g,B}=2$.} \label{fig:SSRvsP_Q_sCSI}
\end{figure}

Figs. \ref{fig:SSRvsP_Q_pCSI} and \ref{fig:SSRvsP_Q_sCSI} consider a similar scenario as Figs. \ref{fig:SEEvsP_Q_pCSI} and \ref{fig:SEEvsP_Q_sCSI}, but report the performance in terms of SSR, for perfect and statistical CSI, respectively. Similar considerations as for the SEE apply. In particular, using $3$ or more quantization bits causes a negligible performance loss. Instead, when $2$ bits are used, a visible, but acceptable gap emerges, while a significant performance degradation is observed when $1$ quantization bit is used. Moreover, also in terms of SSR, the gap is larger in the perfect CSI scenario. Indeed, in this case a visible gap is observed also when $4$ quantization bits are used. However, the gap emerges only for $P_{t,max}\geq 40\,\textrm{dBm}$, which is a higher power value than the typical transmit power values of mobile nodes.

\begin{figure}[!h]
	\centering
	\includegraphics[width=0.5\textwidth]{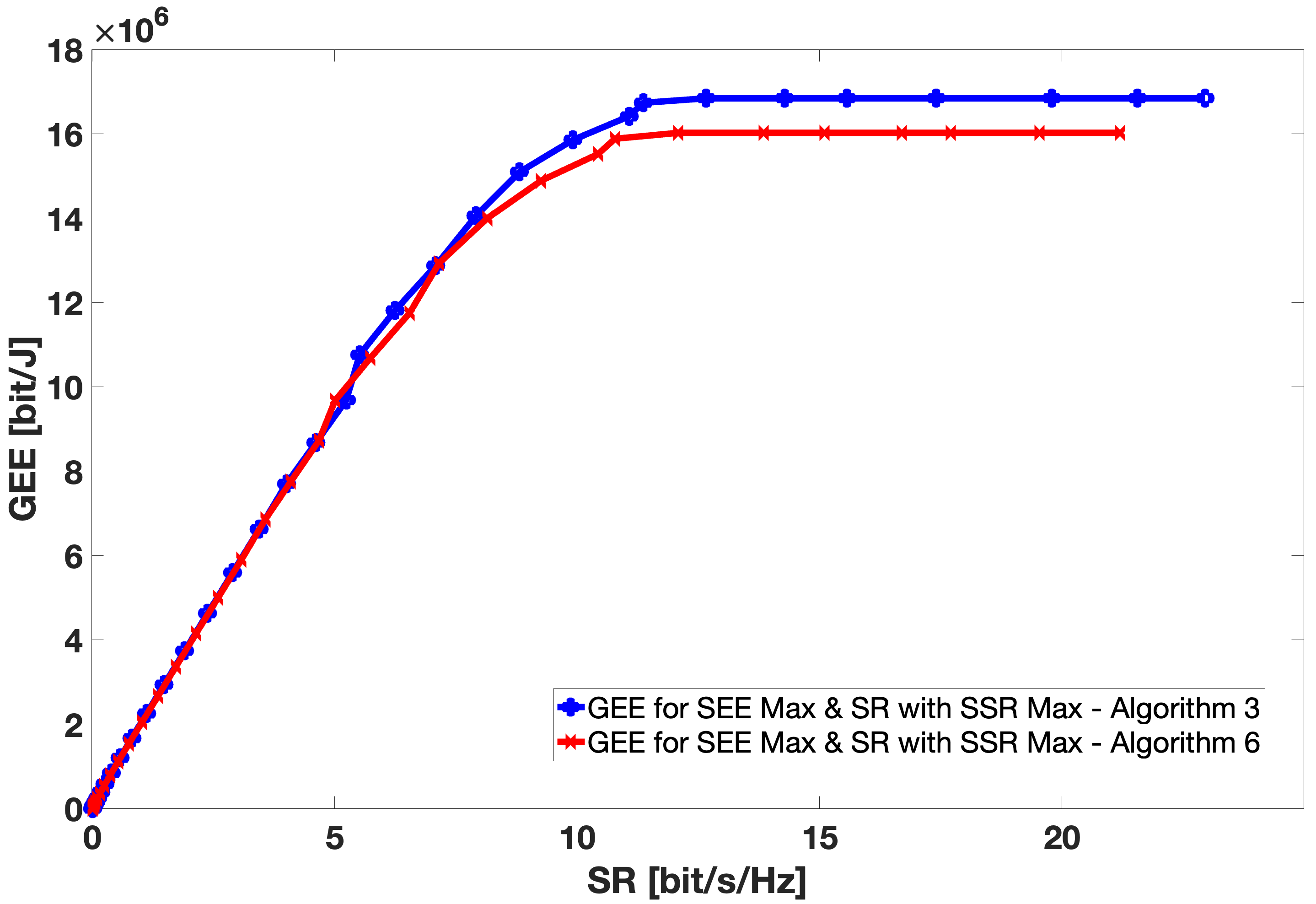}\caption{Achieved GEE for SEE maximization and achieved SR for SSR maximization, with perfect and statistical CSI and $P_{tmax}=[-20, 50]\,\textrm{dBm}$. $K=4$, $N_{B}=4$, $N=100$, $n_{h}=n_{g,E}=4, n_{g,B}=2$.} \label{fig:GEEvsRate}
\end{figure}

Finally, Fig. \ref{fig:GEEvsRate} illustrates the performance of the proposed schemes in terms of global energy efficiency (GEE) and sum-rate (SR). Specifically, each plot point refers to a specific value of the maximum transmit power $P_{tmax}$, for which Algorithms \ref{Alg:SEE1} and \ref{Alg:SEE2} have been run, for both SEE and SSR maximization. Then, the y-axis of Fig. \ref{fig:GEEvsRate} shows the GEE obtained when the system resource are used to maximize the SEE, while the x-axis shows the SR when the system resources are used to maximize the SSR. It is seen that the SR is monotonically increasing, as expected, since SSR maximization leads to using all of the available power $P_{tmax}$. Thus, for increasing $P_{tmax}$, both the SR and SSR will increase. Instead, it is seen that the GEE saturates for a sufficiently large value of $P_{tmax}$. Indeed, both the SEE and GEE are unimodal function in the transmit powers and thus, once the value of $P_{tmax}$ is large enough to attain the peak of the SEE, the resource allocation does not change if $P_{tmax}$ is further increased, and thus both the SEE and GEE saturate.

\section{Conclusion} \label{Sec:Concl} 
This work has studied the problem of SEE maximization in a RIS-aided wireless network. Two provably convergent and low-complexity algorithms have been proposed to allocate the RIS coefficients, users' transmit powers, and receive filters, with both perfect and statistical CSI and the use of both active and nearly-passive RISs. The study has shown that: (1) active RISs do not always yield higher SEE than nearly-passive RISs due to the higher static power consumed by each RIS element; (2) the lack of perfect CSI  does not cause a significant SEE degradation. A promising future line of work is the use of RISs in the near field of network transceivers, also considering electromagnetic compliant models \cite{shabir2025electromagnetically} based on the theory of anomalous reflection \cite{diaz2017generalized}. Another promising line of investigation is the consideration of networks with multiple antennas also at the receiver side.

\bibliographystyle{IEEEtran}
\bibliography{FracProg.bib, references.bib}

\end{document}

%% file: style_icc2024_ext.tex
\usepackage[left=0.625in, right=0.625in, top=0.75in, bottom=1in]{geometry}

% ...

\usepackage{enumitem}
\setlist[itemize]{label=$\bullet$}

\usepackage{hyperref}

\hypersetup{
	colorlinks=true,       
	linkcolor=blue,       
	citecolor=green,      
	filecolor=magenta,    
	urlcolor=cyan         
}

\usepackage{caption}
\usepackage{url}
\makeatletter
\def\url@remove@extrahref#1{%
	\xdef\@savedhref{#1}%
	\expandafter\expandafter\expandafter\strip@url
	\expandafter\@savedhref
	\expandafter\@gobble
}
\makeatother



%% file: acronyms.tex
\makeglossaries

\newacronym{mac}{MAC}{multiple-access channel}
\newacronym{bc}{BC}{broadcast channel}
\newacronym{mimo}{MIMO}{multiple-input multiple-output}
\newacronym{siso}{SISO}{single-input single-output}
\newacronym{sc}{SC}{single-carrier}
\newacronym{mc}{MC}{multi-carrier}
\newacronym{ofdma}{OFDMA}{orthogonal frequency division multiple access}
\newacronym{af}{AF}{amplify-and-forward}
\newacronym{df}{DF}{decode-and-forward}
\newacronym{cf}{CF}{compress-and-forward}
\newacronym{mwrc}{MWRC}{multi-way relay channel}
\newacronym{pde}{PDE}{partial data exchange}
\newacronym{fde}{FDE}{full data exchange}
\newacronym{iid}{i.i.d.\@}{independent and identically distributed}
\newacronym{awgn}{AWGN}{additive white Gaussian noise}
\newacronym{awg}{AWG}{additive white Gaussian}
\newacronym{sic}{SIC}{successive interference cancellation}
\newacronym{snr}{SNR}{signal-to-noise ratio}
\newacronym{sinr}{SINR}{signal-to-interference-plus-noise ratio}
\newacronym{ber}{BER}{bit error rate}
\newacronym{zf}{ZF}{zero-forcing}
\newacronym{mmse}{MMSE}{minimum mean square error}
\newacronym{sud}{SUD}{single user decoding}
\newacronym{dof}{DoF}{degrees of freedom}
\newacronym{gdof}{GDoF}{generalized degrees of freedom}
\newacronym{nnc}{NNC}{noisy network coding}
\newacronym{dmn}{DMN}{discrete memoryless network}
\newacronym{csi}{CSI}{channel state information}
\newacronym{ee}{EE}{energy efficiency}
\newacronym{ian}{IAN}{treating interference as noise}
\newacronym{snd}{SND}{simultaneous non-unique decoding}
\newacronym{brd}{BRD}{best response dynamics}
\newacronym{br}{BR}{best response}
\newacronym{ne}{NE}{Nash equilibrium}
\newacronym{lhs}{LHS}{left-hand side}
\newacronym{rhs}{RHS}{right-hand side}
\newacronym{gee}{GEE}{global energy efficiency}
\newacronym{wsee}{WSEE}{weighted sum energy efficiency}
\newacronym{wpee}{WPEE}{weighted product energy efficiency}
\newacronym{wmee}{WMEE}{weighted minimum energy efficiency}
\newacronym{kkt}{KKT}{Karush-Kuhn-Tucker}
\newacronym{pc}{PC}{pseudo-concave}
\newacronym{qc}{QC}{quasi-concave}
\newacronym{ql}{QL}{quasi-linear}
\newacronym{pl}{PL}{pseudo-linear}
\newacronym{spc}{SPC}{strictly pseudo-concave}
\newacronym{sqc}{SQC}{strictly quasi-concave}
\newacronym{lfp}{LFP}{linear fractional problem}
\newacronym{clfp}{CLFP}{concave-linear fractional problem}
\newacronym{ccfp}{CCFP}{concave-convex fractional problem}
\newacronym{mmfp}{MMFP}{max-min fractional problem}
\newacronym{sorp}{SoRP}{sum-of-ratios problem}
\newacronym{porp}{PoRP}{product-of-ratios problem}
\newacronym{srp}{SRP}{single-ratio problem}
\newacronym{brb}{BRB}{branch-reduce-and-bound}
\newacronym{qos}{QoS}{quality-of-service}
\newacronym{comp}{CoMP}{cooperative multi-point}
\newacronym{ue}{UE}{user equipment}
\newacronym{bs}{BS}{base station}
\newacronym{mrc}{MRC}{maximum ratio combining}
\newacronym{d2d}{D2D}{device-to-device}
\newacronym{lmmse}{LMMSE}{linear minimum mean square error}
\newacronym{ris}{RIS}{reconfigurable intelligent surface}
\newacronym{svd}{SVD}{singular values decomposition}